\documentclass{article}

\usepackage{PRIMEarxiv}

\usepackage[utf8]{inputenc} % allow utf-8 input
\usepackage[T1]{fontenc}    % use 8-bit T1 fonts
\usepackage{hyperref}       % hyperlinks
\usepackage{url}            % simple URL typesetting
\usepackage{booktabs}       % professional-quality tables
\usepackage{amsfonts}       % blackboard math symbols
\usepackage{amsmath}
\usepackage{nicefrac}       % compact symbols for 1/2, etc.
\usepackage{microtype}      % microtypography
\usepackage{lipsum}
\usepackage{graphicx}
\usepackage{subfigure}
\usepackage{mathpazo}
\usepackage{sectsty}
\usepackage{tikz}

\usepackage{booktabs}
\usepackage{subcaption}
\usepackage{multirow}
\usepackage{array}
\usepackage{color}
\newcommand{\change}[1]{\textcolor{black}{#1}}

\usepackage{pifont}% http://ctan.org/pkg/pifont
\newcommand{\cmark}{\ding{51}}%
\newcommand{\xmark}{\ding{55}}%

\graphicspath{{media/}}     % organize your images and other figures under media/ folder

\def\bfd{{\mathbf{d}}}
\def\bff{{\mathbf{f}}}
\def\bfg{{\mathbf{g}}}

\def\bfu{{\mathbf{u}}}
\def\bfv{{\mathbf{v}}}
\def\bfx{{\mathbf{x}}}

\def\bfz{{\mathbf{z}}}
\def\bfdelta{{\boldsymbol{\delta}}}

\def\bfxi{{\boldsymbol{\xi}}}

\def\calG{\mathcal{G}}

\def\calU{\mathcal{U}}

\title{Extending Parametric Model Embedding with Physical Information for Design-space Dimensionality Reduction in Shape Optimization
}

\author{
  A. Serani$^{1,\star}$, G. Palma$^1$, J. Wackers$^2$, D. Quagliarella$^3$, S. Gaggero$^4$, and M. Diez$^1$\\
  $^1$National Research Council-Institute of Marine Engineering, Rome, Italy\\
  $^2$Ecole Centra de Nantes, Nantes, France\\
  $^3$Italian Aerospace Research Centre, Capua, Italy\\
  $^4$University of Genoa, Genoa, Italy\\
  $^\star$\texttt{andrea.serani@cnr.it} \\
}

\begin{document}
% \begin{tikzpicture}[remember picture,overlay]
%    \node [rectangle, fill=cyan, fill opacity=0.5, anchor=north, minimum width=\paperwidth, minimum height=3cm, text width=\textwidth, align=center, text height=5ex, text depth=10ex, align=left] at (current page.north) {\sffamily\small 
%    \textbf{This is a preprint of the following article:}\\
%    A. Serani, T. P. Scholcz, V. Vanzi, A Scoping Review on Simulation-based Design Optimization in Marine Engineering: Trends, Best Practices, and Gaps. \textit{Archives of Computational Methods in Engineering}, 2024.\\
%    \textbf{The published article is available by following the DOI: \texttt{10.1007/s11831-024-10127-1}, which may differ from this preprint.}
%    };
% \end{tikzpicture}

\begin{tikzpicture}[remember picture,overlay]
   % Nodo per il riempimento con trasparenza
   \node [rectangle, fill=cyan, fill opacity=0.5, anchor=north, minimum width=\paperwidth, minimum height=3cm] at (current page.north) {};

   % Nodo separato per il testo, senza trasparenza
   \node [anchor=north, minimum width=\paperwidth, minimum height=3cm, text width=\textwidth, align=center, text height=5ex, text depth=15ex, align=left] at (current page.north) {
     \sffamily\small
     \textbf{This is a preprint submitted to:} \textit{Engineering with Computers}
     % \textbf{This is a preprint of the following article:}\\
     % A. Serani and M. Diez, A Scoping Review on Simulation-based Design Optimization in Marine Engineering: Trends, Best Practices, and Gaps. \textit{Archives of Computational Methods in Engineering}, 2024.\\
     % \textbf{The published article is available by following the DOI: \texttt{10.1007/s11831-024-10127-1}, which may differ from this preprint.}
   };
\end{tikzpicture}

\maketitle

\begin{abstract}
Design-space dimensionality reduction is essential to mitigate the cost of high-fidelity simulation-based optimization, especially when dealing with high-dimensional geometric parameterizations. Traditional linear techniques, such as principal component analysis, are widely used but often neglect the physical response of the system and lack invertibility to the design space, i.e., the ability to reconstruct the original design parameters from a reduced representation.
This work introduces two physics-aware extensions of the parametric model embedding (PME) framework, aimed at generating reduced representations that incorporate physical information while maintaining analytical backmapping. The first, physics-informed PME (PI-PME), combines geometric and physical variability; the second, physics-driven PME (PD-PME), relies solely on physical responses.
The proposed methods enable the construction of interpretable and physically relevant reduced spaces that can be used for design-space exploration, surrogate modeling, and optimization. The approach is demonstrated on multiple engineering configurations, including airfoils, propellers, gliders, and hulls, showing its ability to capture performance-relevant directions and preserve parametric consistency.
The methodology is offline and non-intrusive, compatible with low-fidelity simulations, and requires only a modest number of samples to ensure variance convergence
\end{abstract}

% keywords can be removed

\section{Introduction}
Shape optimization of functional surfaces presents a multifaceted challenge, characterized by numerous geometric, functional, and performance constraints. This is particularly evident in vehicle design, where modern processes must address a broad spectrum of requirements—ranging from energy efficiency and cost effectiveness to safety and environmental sustainability—while simultaneously accounting for aerodynamics, hydrodynamics, structural integrity, noise emissions, and regulatory standards. As the complexity of each vehicle concept increases, so does the dimensionality of the design space, often leading to an exponential growth in the number of parameters that define the shape and operating conditions. This phenomenon, commonly referred to as the \textit{curse of dimensionality} \cite{bellman1957dynamic}, complicates the exploration, analysis, and optimization tasks, since an extremely large solution space quickly becomes prohibitive to sample thoroughly.

Design-space dimensionality reduction methods for shape optimization \cite{serani2024survey} have been extensively explored to address these challenges, simplifying high-dimensional design spaces without significantly compromising predictive accuracy \change{and can be broadly categorized into \textit{offline} and \textit{online} approaches. Offline or upfront methods (e.g., principal component analysis and its variants \cite{yonekura2014shape,diez2015-CMAME}) aim to construct a reduced representation before the actual optimization or design exploration process, relying on a fixed dataset of samples. Conversely, online or adaptive techniques (e.g., active subspace methods \cite{constantine2014active,berguin2015dimensionality}) seek to refine the reduced space on-the-fly, as new data becomes available. While online approaches may offer enhanced flexibility, they often require substantial computational overhead and complex integration with the optimization workflow. In contrast, offline methods are appealing for their simplicity, interpretability, and suitability for use with high-fidelity simulations, provided that the initial sampling ensures sufficient coverage of the design space.
Among them,} linear methods, such as principal component analysis (PCA) \cite{yonekura2014shape,harries2021application,jun2020application,zhang2024geometric} and singular value decomposition (SVD) \cite{poole2017high,allen2018wing,poole2022efficient} derived from proper orthogonal decomposition (POD) \cite{zhang2018multidisciplinary,ballarin2019pod,yanhui2019performance,yamazaki2020efficient,yang2024aerodynamic,li2024aerodynamic} and Karhunen-Loève expansion (KLE) \cite{diez2015-CMAME,d2020design,chang2023research}, have been widely adopted for reducing the dimensionality of design space by capturing the dominant modes of geometric variance. However, traditional PCA-based approaches primarily utilize geometric data, potentially neglecting critical physical phenomena influencing design performance. To enhance predictive capabilities, physics-informed PCA methods have recently emerged \cite{serani2022hull,diez2023design,li2023new,zhang2023efficient,mingzhi2024breaking}, augmenting geometric data matrices with physical simulation outputs, such as lumped parameters (e.g., efficiency or drag coefficients) or distributed physical quantities (e.g., pressure or velocity fields). Furthermore, \cite{khan2022shape, khan2022geometric, masood2023shape, kostas2023machine,masood2024generative} introduced geometric moments as cost-effective, physics-related descriptors that substitute direct simulation data, further enriching the PCA model without the computational cost of numerical simulations. These enriched representations allow dimensionality reduction frameworks to capture not only geometric but also essential physical variability, significantly improving design relevance and effectiveness.

\change{Despite their widespread use, conventional offline dimensionality reduction techniques such as PCA and SVD present critical limitations in design applications. First, SVD-based approaches, which are often employed in POD and KLE formulations, do not support weighting of the inputs due to their inherent algebraic structure. This prevents users from tuning the relative importance of geometric features or physical observables during the reduction.
In contrast, PCA—formulated as an eigenvalue problem of a covariance matrix—can accommodate weighted formulations, enabling generalized PCA variants and physics-informed extensions that incorporate additional information beyond shape geometry. However, these approaches still lack a crucial capability: the ability to reconstruct a feasible set of design parameters from a point in the reduced space. This issue, commonly referred to as the \emph{pre-image problem} \cite{gaudrie2020modeling}, severely limits their applicability in optimization workflows or CAD-integrated environments where direct parametric interpretation and reversibility are required.}

\change{To overcome the latter limitation, the parametric model embedding (PME) technique \cite{serani2023parametric} was introduced. PME embeds the original parametric variables explicitly into the reduction process and re-parameterizes \cite{serani2024aerodynamic} the design space by constructing an orthonormal basis that supports both weighted formulations and analytical backmapping, thereby bridging the gap between statistical reduction and design interpretability.}

\begin{table*}[!t]
    \centering
    \change{ \caption{\change{Comparison of dimensionality reduction methods in terms of data types handled, use of design variables, availability of weight control, and possibility of backmapping to the original parametric space.}}
    \label{tab:comparison_methods}
    \begin{tabular}{lcccccp{4.25cm}}
\toprule
Method & Geometry & Physics & Variables & Weights & Backmapping & Notes \\ %\raggedright Notes \\
         \midrule
    SVD        & \cmark   & --      & --        & \xmark& \xmark  & Standard decomposition of geometry-only snapshots. \\
    PCA        & \cmark   & --      & --        & \cmark& \xmark  & Can weight geometric features. \\
    PI-PCA     & \cmark   & \cmark  & --        & \cmark& \xmark  & Can also incorporate physics thanks to weights. \\
    \midrule
    PME        & \cmark   & --      & \cmark    & \cmark& \cmark  & Embeds geometry and  variables. \\
    PI-PME     & \cmark   & \cmark  & \cmark    & \cmark& \cmark  & Combines geometry and physics with variables embedding. \\
    PD-PME     & --       & \cmark  & \cmark    & \cmark& \cmark  & Fully physics-driven embedding. \\
    \bottomrule
    \end{tabular}
    }
\end{table*}
\change{Nevertheless, the original PME was designed to address the pre-image problem only, without accounting for physical integral quantities or distributed fields, which are essential to characterize performance in fluid-dynamic applications. The present work extends PME to incorporate such physical information, exploring two complementary linear offline strategies: physics-informed PME (PI-PME), which embeds physical outputs alongside geometry and design variables, and physics-driven PME (PD-PME), which relies solely on physical features to define the embedding.}

PI-PME enriches the data matrix with physical information, including lumped parameters (e.g., force coefficients or efficiencies) or distributed quantities (e.g., surface pressure fields), derived from low-fidelity simulations. Although this enrichment entails an upfront computational cost, the resulting latent space captures both geometric and physical variability, enabling a more informed and reliable assessment of performance trade-offs. By integrating physical insights at the earliest stages of design exploration, PI-PME provides a foundation that can guide subsequent high-fidelity simulations and optimizations.

Building on the same principle, a PD-PME is also proposed, which goes one step further by focusing solely on physical parameters and effectively excluding geometric variance. Although this approach removes explicit shape information from the dimensionality reduction, it can be advantageous in scenarios where the primary objective is to characterize or optimize physical behaviors, and where the link between geometry and physics is well captured through simplified models. This exclusive emphasis on physics can also serve as a natural starting point for multi-fidelity optimization or for data-driven reduced-order modeling, where lumped or distributed physical data are used to predict key performance indicators at reduced computational expense.

\change{Table~\ref{tab:comparison_methods} summarizes the distinguishing features of the considered dimensionality reduction methods with respect to four key aspects: (i) the type of data handled (geometry, physics, or both), (ii) the inclusion of design variables in the embedding, (iii) the availability of weight control, and (iv) the capability to perform backmapping to the original parametric space. Classical methods such as SVD and PCA rely solely on geometric information and, while PCA can incorporate weights to prioritize specific features, neither provides a mechanism for recovering design parameters from reduced coordinates. Physics-informed variants of PCA extend this formulation by integrating physical observables—either lumped or distributed—but still lack the ability to reconstruct input designs.
In contrast, the PME framework embeds design variables directly into the reduction process, enabling analytical backmapping and supporting weighted formulations. Its extensions, PI-PME and PD-PME, further enrich the embedding by incorporating physical quantities. While PI-PME retains geometric and parametric information, PD-PME constructs the latent space purely from physical descriptors, which can be advantageous when physical trends are more relevant than shape variance. These distinctions highlight the versatility of PME-based approaches in design-space reduction tasks where physical interpretability and parameter reconstruction are essential.}

The effectiveness of PI-PME and PD-PME becomes particularly apparent in real-world applications such as bio-inspired underwater gliders, naval propellers, or classical airfoils (e.g., the RAE-2822). In these cases, fluid-dynamic effects strongly impact performance, and omitting them from the dimensionality reduction step may lead to suboptimal designs and increased resource expenditure. By embedding physical information—whether partially, as in PI-PME, or exclusively, as in PD-PME—into the reduced space, it becomes possible to filter the vast design domain more intelligently, ensuring that only the most relevant configurations are explored. This synergy between geometric and physical factors not only improves the quality of early design-stage decisions but also streamlines the path to final optimization, cutting down both computational time and overall development costs. 
\change{It is important to emphasize that the test cases discussed throughout the paper are not simplified academic benchmarks, but rather representative of real-world early-stage design workflows in aerospace and marine engineering. Each case reflects a high-dimensional design problem where physical performance plays a central role: the RAE-2822 airfoil is a canonical transonic shape studied for aerodynamic efficiency; the seakeeping hull and propeller designs involve complex hydrodynamic effects under realistic operating conditions; and the glider case addresses the viscous performance of bioinspired geometries. In all scenarios, dimensionality reduction is not an optional step, but a necessary enabler to manage complexity, improve interpretability, and make downstream optimization tractable.}
PME’s flexible structure also makes it directly applicable to structural optimization. Instead of relying on aerodynamic or hydrodynamic simulations, one may incorporate structural properties (e.g., mass moments of inertia, stiffness distributions, or stress fields) into the PME framework, enabling integrated assessments of mechanical performance \cite{diez2023design,pellegrini2025-411}.

The following sections illustrate the mathematical underpinnings of PME, PI-PME, and PD-PME, including detailed discussions of how generalized PCA is adapted to incorporate physical data, and how these frameworks can be applied to a range of high-dimensional vehicle design problems. Through examples, the benefits of early-stage integration of physical insights will be demonstrated, showing improved design robustness, a decrease in exhaustive computational campaigns, and the potential for more agile innovation cycles in modern vehicle development. It may finally be noted that, while this work focuses on vehicle applications, the proposed method is broadly applicable to general shape optimization problems involving functional surfaces influenced by physical constraints.

\section{Design-space dimensionality reduction methods}
In shape optimization problems aimed at minimizing an objective function $f(\bfu)$, the optimum shape is uncertain before and during the optimization process. The design variable vector $\bfu \in \mathcal{U} \subset \mathbb{R}^M$, with $\mathcal{U}$ the design domain, can therefore be treated as a random variable, with an associated probability density $p(\bfu)$ which reflects the likelihood of a certain design being optimal, given the existing knowledge. If nothing is known about the design, $p$ can be uniform. This uncertainty propagates through the modeling process, generating uncertainties in both geometric configurations and the associated physical responses.

\change{Dimensionality reduction techniques such as PCA, SVD, POD, and KLE provide a systematic way to extract the dominant directions of variability from such data, i.e., the combinations of design variables which most influence the results. These methods enable the construction of compact, low-dimensional representations that preserve the most informative features for design analysis and decision-making.}

As will be discussed in the following sections, integrating physical insights into these dimensionality reduction frameworks, for design-space assessment and subsequent shape optimization, provides a more robust path toward identifying optimal solutions under uncertainty.

\subsection{Parametric model embedding}
{PME \cite{serani2023parametric} is a design-space dimensionality reduction method that extends the standard PCA approach by incorporating both shape deformations and design variables into a generalized feature space. Specifically, PME \change{solves a generalized eigenvalue problem, applying a generalized} PCA to an augmented matrix that includes the discretized shape deformation vector $\bfd$ and the original design variables $\bfu$. This extension allows PME to directly map the reduced design space back to the original design variables without the need for reparameterization, which is typically required in standard PCA approaches. As a result, PME offers a more robust and practical method for maintaining the integrity of the original design features while facilitating effective shape optimization.} 

\begin{figure}[!b]
    \centering
    \includegraphics[width=0.5\columnwidth]{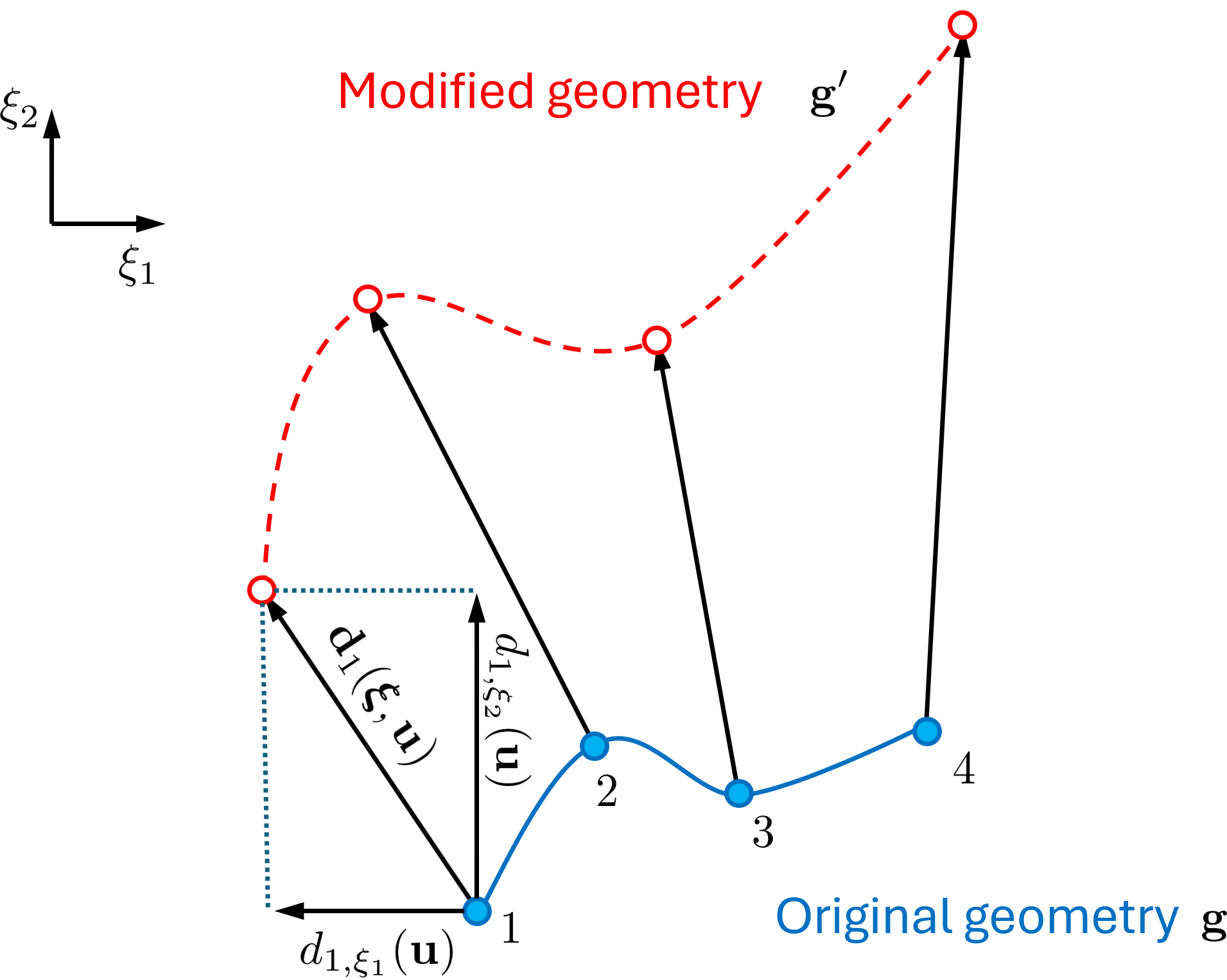}
    \caption{Example of shape modification and discretization with notation, where $n=2$ and $L=4$}
    \label{fig:shape_ex}
\end{figure}
Consider a manifold $\calG$, which identifies the original shape/geometry. {This manifold defines the geometric space in which the shape is parameterized by curvilinear coordinates $\bfxi\in\calG$. The} coordinates {of the original shape} are represented by $\bfg(\bfxi)\in\mathbb{R}^n$ with $n=1,2,$ or 3. Assume that, for~the purpose of shape optimization, $\bfg$ can be transformed to a deformed shape/geometry $\bfg'(\bfxi,\bfu)$ by
\begin{equation}
\bfg'(\bfxi,\bfu)=\bfg(\bfxi)+\bfdelta(\bfxi,\bfu) \qquad \forall \bfxi\in\calG
\end{equation} 
where $\bfdelta(\bfxi,\bfu)\in\mathbb{R}^n$ is the resulting shape modification vector, defined by arbitrary shape parameterization or modification method (e.g., CAD parameterization, Bezier surfaces, free-form deformation, NURBS, etc.). %, and $\bfu$\in\calU\subset\mathbb{R}^M$  is the design-variable vector of dimension $M$. 

Discretizing $\calG$ by $L$ elements (see, e.g., Fig. \ref{fig:shape_ex}) of measure $\Delta\calG_j$ (with $j=1,\dots,L$), and sampling $\calU$ by a statistically convergent number of realizations $S$, so that $\{\bfu_k\}_{k=1}^S  \sim p(\bfu)$, ${\mathbf{d}}(\bfxi,\bfu)$ can be obtained as the discretization of ${\bfdelta}(\bfxi,\bfu)$. Organizing $\hat{\mathbf{d}}=\bfd - \langle\bfd\rangle$ (with $\langle \cdot \rangle$ the mean value) in a data matrix $\mathbf{D}$ of dimensionality $\left[nL\times S\right]$, one obtains
\begin{equation}\label{eq:data}
\begin{aligned}
\mathbf{D} &=
\left[
\begin{array}{ccc}
\hat{\bfd}_{1}^{(1)} & \dots & \hat{\bfd}_{1}^{(S)}\\
\vdots & \vdots & \vdots\\
\hat{\bfd}_{n}^{(1)} & \dots & \hat{\bfd}_{n}^{(S)}\\
\end{array}
\right]\\
& = 
\left[
\begin{array}{ccc}
\hat{d}_{1,\xi_1}(\bfu_1) & & \hat{d}_{1,\xi_1}(\bfu_S)\\
\vdots & & \vdots\\
\hat{d}_{L,\xi_1}(\bfu_1) & & \hat{d}_{L,\xi_1}(\bfu_S)\\
\vdots & \dots & \vdots\\
\hat{d}_{1,\xi_n}(\bfu_1) & & \hat{d}_{1,\xi_n}(\bfu_S)\\
\vdots & & \vdots\\
\hat{d}_{L,\xi_n}(\bfu_1) & & \hat{d}_{L,\xi_n}(\bfu_S)\\
\end{array}
\right]
\end{aligned}
\end{equation}   
where $\hat{d}_{j,\xi_k}$ is the $k$-th component of the shape modification vector associated to the $j$-th element.
The~embedding is achieved by defining the matrix $\mathbf{P}$ of dimensionality $\left[(nL+M)\times S\right]$ as follows
\begin{equation}\label{eq:P}
\mathbf{P}=\left[
\begin{array}{c}
\mathbf{D}\\
\mathbf{U}
\end{array}
\right]
%\qquad
\end{equation}
with
\begin{equation}
\mathbf{U}=
\left[
\begin{array}{ccc}
\hat{\bfu}^{(1)} & \dots & \hat{\bfu}^{(1)}\\ 
\end{array}
\right]
=
\left[
\begin{array}{ccc}
\hat{u}_{1,1} & & \hat{u}_{1,S}\\
\vdots & \cdots & \vdots\\
\hat{u}_{M,1} & & \hat{u}_{M,S}\\
\end{array}
\right]
\end{equation}   
where $\hat{\bfu}=\bfu-\langle\bfu\rangle$.
The matrix $\mathbf{U}$ is appended to the data matrix $\mathbf{D}$ and associated to a null weight {$\mathbf{W}_\bfu$ such that
\begin{equation}\label{eq:pme_weight}
\mathbf{W}_\bfu=\mathbf{0}
\qquad
\mathrm{and}
\qquad
\mathbf{W}=\left[
\begin{array}{cc}
\mathbf{W}_\bfd & \mathbf{0}\\
\mathbf{0} & \mathbf{W}_\bfu\\
\end{array}
\right]
\end{equation}   
}
and so leading to a generalized PCA problem in the form
\begin{equation}\label{eq:pme_pca}
\mathbf{A}\mathbf{G}\mathbf{W}\mathbf{Z}=\mathbf{Z}\boldsymbol{\Lambda} \qquad \mathrm{with} \qquad {\mathbf{A}}=\frac{1}{S}\mathbf{PP}^\mathsf{T} 
\end{equation} 
where
\begin{equation}\label{eq:pme}
\mathbf{G}=\left[
\begin{array}{cc}
\mathbf{G}_\bfd & \mathbf{0}\\
\mathbf{0} & \mathbf{I}\\
\end{array}
\right]
\end{equation}
and
\begin{equation}
{\mathbf{Z}} = \left[{\bfz}_1 \,\,\, \dots \,\,\, {\bfz}_{S} \right]
%\end{equation}
\qquad
\mathrm{with}
\qquad
%\begin{equation}
{\bfz}_k=\left[
\begin{array}{c}
\mathbf{q}_k\\
\bfv_k
\end{array}
\right]
\end{equation} 
\change{is the eigenvectors matrix, whereas $\boldsymbol{\Lambda}$ is the associated eigenvalues matrix.}
{{Here}, $\mathbf{q}_k$ and $\bfv_k$ represent the eigenvector components associated to the shape modification $\bfd$ and design variable $\bfu$ vectors, respectively.}
The matrix $\mathbf{G}_\bfd=\mathrm{diag}\left(\mathbf{G}_1,\dots,\mathbf{G}_n \right)$ is block diagonal and has dimensionality $\left[nL\times nL\right]$, with~each $\left[L\times L\right]$ $k$-th block being a diagonal matrix itself
\begin{equation}
\mathbf{G}_k=\mathrm{diag}\left(\Delta\calG_1, \dots, \Delta\calG_L\right) 
\end{equation}
containing the measure $\Delta\calG_j$ of the $j$-th element. Similarly, 
$\mathbf{W}_\bfd=\mathrm{diag}\left(\mathbf{W}_1,\dots,\mathbf{W}_n \right)$ is a block diagonal matrix of dimensionality $\left[nL\times nL\right]$, where each $\left[L\times L\right]$ $k$-th block $\mathbf{W}_k$ ($k=1,\dots,n$) is itself a diagonal matrix defined as
\begin{equation}\label{eq:weights}
\mathbf{W}_k=\mathrm{diag}\left({\rho_1}, \dots, {\rho_L}\right)
\end{equation}
where $\rho_j$  (for $j = 1,\dots, L$) represents the arbitrary weight given to each~element.
The columns of $\mathbf{Z}$ are normalized to unit norm with respect to the $\mathbf{G}\mathbf{W}$ scalar product. Specifically, each column $\mathbf{z}_k$ is scaled by a scalar $\gamma_k$ such that:
\begin{equation}
\gamma_k = \sqrt{\mathbf{z}_k^\mathsf{T} \mathbf{G} \mathbf{W} \mathbf{z}_k}, \qquad \mathbf{z}_k^\star = \frac{\mathbf{z}_k}{\gamma_k},
\end{equation}
leading to the normalized matrix:
\begin{equation}
\mathbf{Z}^\star = \mathbf{Z} \mathbf{\Gamma}^{-1},
\end{equation}
where
\begin{equation}
\mathbf{\Gamma} = \mathrm{diag}\left( \left[ \sqrt{ \mathbf{z}_1^\mathsf{T} \mathbf{G} \mathbf{W} \mathbf{z}_1 }, \dots, \sqrt{ \mathbf{z}_{nL}^\mathsf{T} \mathbf{G} \mathbf{W} \mathbf{z}_{nL} } \right] \right).
\end{equation}
This normalization ensures that each eigenvector contributes equally, avoiding numerical ill-conditioning due to significant variations in vector norms.

The solutions $\lambda_k$ and the corresponding normalized eigenvectors $\mathbf{v}_k^\star$ (columns component of $\mathbf{Z}^\star$ associated to the original design variables, see Eq. \ref{eq:pme}) are used to construct the reduced dimensionality representation of the original parameterization by means of the $N$ reduced design variables $\mathbf{x} = [x_1 , \dots , x_N]^\mathsf{T}$. Defining the desired confidence level $l$, with~$0<l\leq1$, the~number of reduced design variables $N$ is chosen such that
\begin{equation}\label{eq:lambda}
\sum_{k=1}^N \lambda_k\geq l\sum_{k=1}^{nL} \lambda_k=l\sigma^2 \quad \mathrm{with} \quad \lambda_k\geq\lambda_{k+1},
\end{equation}
and the PME of the original design variables is achieved by using these normalized eigenvectors as follows
\begin{equation}\label{eq:recu}
\mathbf{u} \approx \check{\mathbf{u}}=\langle{\mathbf{u}}\rangle+\sum_{k=1}^N x_k\mathbf{v}_k^\star.
\end{equation}
\change{
Thus design parameters $\mathbf{u}$ are assigned zero weight in the inner product definition, but they remain structurally embedded in the generalized eigenvalue problem. This allows PME to retain a direct mapping from the reduced space back to the original parametric domain, thus enabling analytical backmapping (see Eq. \ref{eq:recu}). As a result, the dimension reduction does not compromise the ability to reconstruct feasible and interpretable design configurations.
}

\change{Finally,} to ensure that all the samples in $\mathbf{P}$ can be reconstructed through the reduced-order representation of the original design space, the reduced design variables $\bfx$ can be bounded by identifying the minimum and maximum values attained by each component $\Theta_{jk}$ of the projection coefficients $\boldsymbol{\theta}_j$, as follows: 
\begin{equation}\label{eq:theta} 
\min_j \Theta_{jk} \leq x_k \leq \max_j \Theta_{jk} \qquad \text{for } k = 1, \dots, N, 
\end{equation} 
with $\boldsymbol{\Theta}=[\boldsymbol{\theta}_1 , \dots , \boldsymbol{\theta}_S]^\mathsf{T}$ evaluated by projecting the matrix $\mathbf{P}$ onto the normalized basis $\mathbf{Z}^\star$, i.e., 
\begin{equation} \boldsymbol{\Theta} = \mathbf{P}^\mathsf{T} \mathbf{G} \mathbf{W} \mathbf{Z}^\star. 
\end{equation}
Moreover, it can be shown that the sum of the squared projection coefficients across the dataset equals the sum of the corresponding eigenvalues \cite{diez2015-CMAME}: 
\begin{equation} 
\dfrac{1}{S}\sum_{j=1}^S \Theta_{jk}^2 = \lambda_k \qquad \text{for } k = 1, \dots, N. \end{equation}
This implies that the variance captured along each principal direction is preserved in the projection. 

{It} may be noted that the overall methodology is independent of the specific shape modification method, which is seen as a black box by PME. 

\subsection{Physics-informed parametric model embedding}

Introducing physical information within PME formulations translates into the definition of, similarly to Eq. \ref{eq:P}, a new matrix $\mathbf{P}_I$ as
\begin{equation}
\mathbf{P}_I=\left[
\begin{array}{c}
\mathbf{D}\\
\mathbf{U}\\
\mathbf{F}\\
\mathbf{C}
\end{array}
\right]
\end{equation}
with
\begin{equation}
\mathbf{F}=
\left[
\begin{array}{ccc}
\vert & \vdots & \vert\\
\hat{\bff}_{j}^{(1)} & \dots & \hat{\bff}_{j}^{(S)}\\
\vert & \vdots & \vert\\
\end{array}
\right]
\end{equation}
and
\begin{equation}
\mathbf{C}=
\left[
\begin{array}{ccc}
\vert & \vdots & \vert\\
\hat{c}_{j}^{(1)} & \dots & \hat{c}_{j}^{(S)}\\
\vert & \vdots & \vert\\
\end{array}
\right]
\end{equation}
where $\mathbf{F}$ is a matrix containing distributed physical information (with $\hat{\bff}_j=\bff_j-\langle\bff_j\rangle$), such as, e.g., pressure distribution and/or velocity components, that don't necessarily have to be defined on the geometry surface but can also belong to the field surrounding the object, like, e.g., wake and vortices, whereas $\mathbf{C}$ contains lumped or scalar quantities of interest (with $\hat{c}_j=c_j-\langle c_j\rangle$), such as, e.g., lift and drag forces.
It may be noted that for a given geometry, the physical information can be collected for more than one operating condition.

Analogously, a corresponding block-diagonal weight matrix ${\mathbf{W}}_I$ is introduced: 
\begin{equation}\label{eq:pipme_weight}
\begin{aligned}
{\mathbf{W}}_I &=\left[
\begin{array}{cccc}
\mathbf{W}_\bfd & \mathbf{0}      & \mathbf{0}& \mathbf{0}\\
\mathbf{0} & \mathbf{W}_\bfu & \mathbf{0}& \mathbf{0}\\
\mathbf{0} & \mathbf{0} & \mathbf{W}_\bff & \mathbf{0} \\
\mathbf{0} & \mathbf{0} & \mathbf{0}& \mathbf{W}_c\\
\end{array}
\right] \\
& = 
\mathrm{diag}
\Bigl(
\mathbf{W}_\bfd,
\underbrace{\mathbf{0}}_{\mathbf{u}},
\mathbf{W}_{\mathbf{f}},
\mathbf{W}_{c}
\Bigr),
\end{aligned}
\end{equation} 
where $\mathbf{W}_\bfd$ is the block-diagonal matrix weighting the geometric entries (as in standard PME), $\mathbf{W}_\bfu$ the one for design variables, $\mathbf{W}_\bff$ for distributed physical data $\mathbf{F}$, and $\mathbf{W}_c$ for lumped scalars $\mathbf{C}$. 

%\change{Each data block in the embedding (geometry, physical fields, or lumped quantities) is normalized via the corresponding weight matrix (see Eq. \ref{eq:pipme_weight}), constructed by inverting the sample variance of its components. This ensures that all quantities—regardless of physical units or magnitude—contribute comparably to the embedding.} 
%
Each block has to be set up to normalize its respective data by the inverse of the estimated variance. 
In practice, for a block associated with a data vector $y\in\{\bfd,\bff,c\}$, one calculates
\begin{equation}\label{eq:var}
\mathrm{Var}(y)=\frac{1}{S}\sum_{j=1}^S (y_j-\langle y \rangle)^2,    
\end{equation}
and then defines
\begin{equation}\label{eq:var_W}
\rho_y = 1/\mathrm{Var}(y).    
\end{equation}
If the block contains multiple components (e.g., multi-dimensional fields or multiple operating conditions), each row or column of that block can be assigned its own weight $\rho_i$. %Doing so ensures that geometric modifications and physical data all contribute comparably to the subsequent PCA procedure—i.e., no portion of the data dominates simply because it exhibits a larger raw variance.
\change{This normalization strategy ensures that both geometric modifications and physical observables contribute on a comparable scale to the PCA procedure—preventing any data block from dominating the embedding solely due to differences in units, magnitude, or raw variance.
}

Once the augmented data matrix $\mathbf{P}_I$ and the block-diagonal weighting matrix ${\mathbf{W}}_I$ are assembled, the eigenvalue problem remains analogous to the original PME formulation:
\begin{equation}
{\mathbf{A}_I}{\mathbf{G}_I}{\mathbf{W}_I}{\mathbf{Z}}_I = {\mathbf{Z}}_I \boldsymbol{\Lambda}_I  
\end{equation}
where
\begin{equation}
\mathbf{A}_I = \frac{1}{S} \mathbf{P}_I\mathbf{P}_I^{\mathsf{T}}, 
\qquad
\mathbf{G}_I = 
\left[
\begin{array}{cccc}
\mathbf{G}_\bfd   & \mathbf{0} & \mathbf{0} & \mathbf{0} \\
\mathbf{0}   & \mathbf{I} & \mathbf{0} & \mathbf{0} \\
\mathbf{0}   & \mathbf{0} & \mathbf{G}_\bff& \mathbf{0} \\
\mathbf{0}   & \mathbf{0} & \mathbf{0} & \mathbf{I}\\ 
\end{array}
\right],
\end{equation}
and
\begin{equation}
\mathbf{Z}_I = \left[\bfz_{I,1} \dots \bfz_{I,S} \right],
\qquad
\bfz_{I,k} = 
\left[
\begin{array}{c}
  \mathbf{q}_k     \\
  \bfv_k \\
  \boldsymbol{\phi}_k\\
  \boldsymbol{\pi}_k 
\end{array} \right]_I
\end{equation}

Here, $\mathbf{G}_\bff$ accounts for the element size of the distributed physical vector, which does not necessarily correspond to $\mathbf{G}_\bfd$, while the identity block is also applied to the lumped scalars. The eigenvector solution $\mathbf{Z}_I$ then provides the reduced representation for all data (where $\boldsymbol{\phi}_k$ and $\boldsymbol{\pi}_k$ are the eigenvector components that embed the distributed and lumped physical parameters), ensuring that the new PI-PME basis incorporates both shape and physics with consistent normalization. This yields a lower-dimensional yet physics-enriched space suitable for subsequent design-space exploration, optimization, or multi-fidelity modeling.

\subsection{Physics-driven parametric model embedding}
In some design scenarios, large geometric changes may not necessarily translate into significant variations in the associated physical phenomena. Consequently, the standard geometry-centric approach in dimensionality reduction can yield subspaces that capture a high geometric variance but do not correlate with improved physical performance. PD-PME addresses this issue by removing geometric deformations (\(\mathbf{D}\)) altogether, focusing only on the physical data. Such a formulation aims to isolate and amplify the directions in the design space that have the strongest impact on physically relevant quantities, rather than those with merely large geometric variability.

By excluding \(\mathbf{D}\), the augmented data matrix reduces to 
\begin{equation*}
\mathbf{P}_D
\;=\;
\begin{bmatrix}
\mathbf{U}\\
\mathbf{F}\\
\mathbf{C}
\end{bmatrix}.
\end{equation*}
Because the objective is to let the physical variability guide the principal components, one assembles a block-diagonal matrix
\begin{equation*}
{\mathbf{W}}_D
\;=\;
\mathrm{diag}
\Bigl(
\underbrace{\mathbf{0}}_{\mathbf{u}},
\,\mathbf{W}_\mathbf{f},
\,\mathbf{W}_c
\Bigr),
\end{equation*}
which assigns, as per standard PME and PI-PME, zero weight to the rows corresponding to \(\mathbf{U}\), while normalizing the distributed physical fields \(\mathbf{F}\) and the lumped scalars \(\mathbf{C}\) by the inverse of their respective variances, as per PI-PME (see, eqs. \ref{eq:var} and \ref{eq:var_W}).  

The solution of the generalized PCA problem
\begin{equation}
{\mathbf{A}}_D
\;{\mathbf{G}}_D
\;{\mathbf{W}}_D
\;{\mathbf{Z}}_D
\;=\;
{\mathbf{Z}}_D
\;\boldsymbol{\Lambda}_D,
\end{equation}
where
\begin{equation}
{\mathbf{A}}_D 
\;=\;
\frac{1}{S}\,\mathbf{P}_D\,\mathbf{P}_D^\mathsf{T},
\end{equation}
yields principal components strictly {driven by physics}, since the geometry has been removed and design variables are weighted to zero.  

A relevant consideration arises when using {only} lumped (scalar) physical data. If the number of such scalars is smaller than the number of design variables, the rank of \(\mathbf{P}_D\) will be limited by the fewer of the two, resulting in a forced dimensionality reduction dictated by the size of the scalar data block rather than by a variance-based choice. Consequently, a minimum requirement for PD-PME to {fully} capture variability in the design vector is that the dimensionality of the physical information (e.g., the number of scalar data points) be at least equal to or larger than the number of design variables. Otherwise, the principal components will be constrained by the lower-dimensional block and may fail to represent all relevant modes of variation.

\section{Test cases}
To highlight the effectiveness and flexibility of the proposed approaches (PME, PI-PME, and PD-PME), four different test cases are considered. Each case targets a distinct geometry and set of operating conditions, thus allowing for a comprehensive assessment of how geometric and physical information can be combined (or, in the case of the physics-driven variant, used exclusively) to enable efficient design-space dimensionality reduction. The following subsections detail each test case in turn, describing both the parameterization strategies and the sources of physical data employed for PI-PME and PD-PME.
\begin{figure*}[!t]
    \centering
    \includegraphics[width=0.49\linewidth]{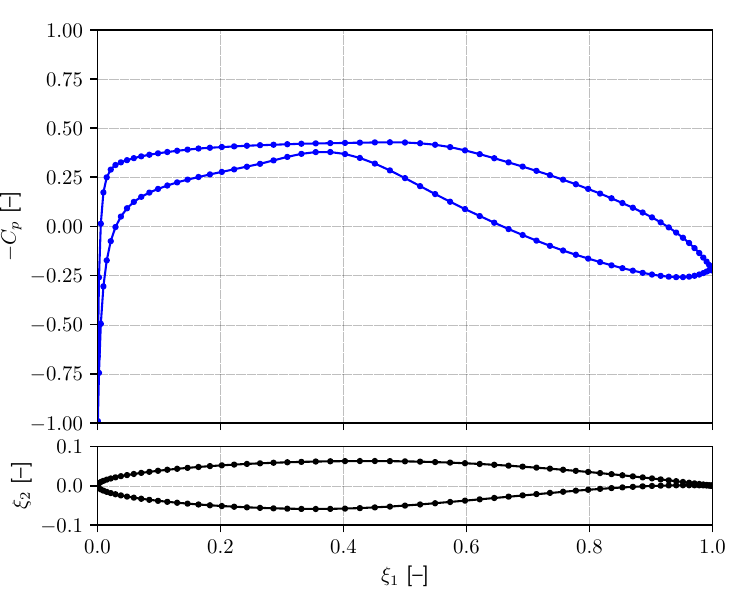}
    \includegraphics[width=0.49\linewidth]{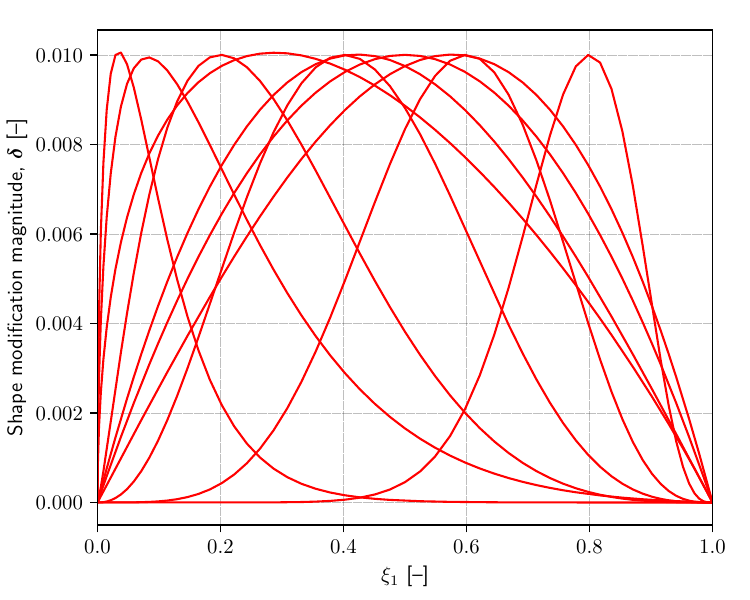}
    \caption{Test case~1: RAE-2822 (left) original geometry and pressure coefficient $C_p$ with discretization and (right) shape modification functions/original parameterization}
    \label{fig:rae_shapes}
\end{figure*}
\begin{figure*}[!b]
    \centering
    \includegraphics[width=1\linewidth]{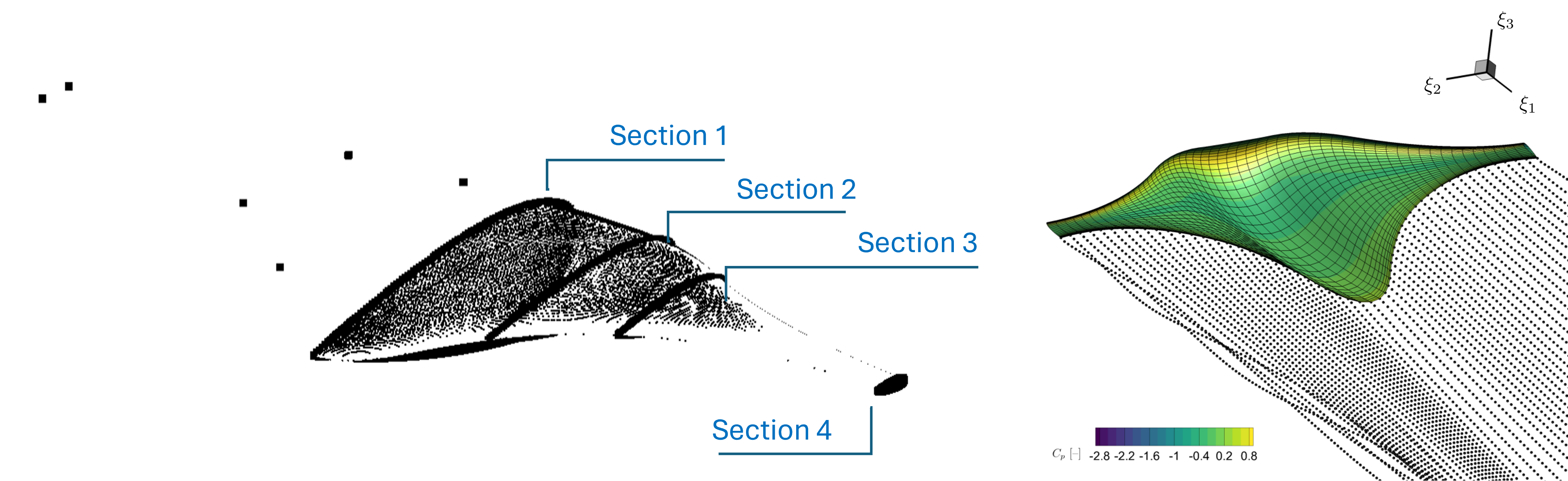}
    \caption{Test case~2: AUG (left) shape parameterization and (right) pressure coefficient and wake as solve outputs}
    \label{fig:glider_shapes}
\end{figure*}

\subsection{Test case 1: RAE-2822 airfoil}
The first test case concerns the design optimization of the classical RAE-2822 airfoil (see Fig.~\ref{fig:rae_shapes}, bottom left). The design space, defined within the activities of the NATO-AVT-331 Research Task Group~\cite{beran2020comparison}, includes \(M = 20\) design variables; each variable is associated with a different shape function (3 polynomials, 6 Hicks--Henne bumps, and 1 Wagner function~\cite{Hicks:78}, see Fig.~\ref{fig:rae_shapes}, right), acting either on the upper or lower surface of the airfoil~\cite{quagliarella202-354}. An in-house code (\texttt{WG2AER}, developed at CIRA) parameterizes the airfoil as a linear combination of the parent geometry \(\bfg(\bfxi)\) and the modification functions \(\bfdelta\).

Discretizing the airfoil with \(L = 129\) grid points, the physical information required to train PI-PME and PD-PME is gathered from the \texttt{XFOIL} solver~\cite{XFOIL}, under operating conditions of Mach~\(=0.4\), Reynolds~\(\,=6.5\times10^6\), and a zero-degree angle of attack. Collected physical data include the pressure coefficient (as distributed information, see Fig.~\ref{fig:rae_shapes} top left) and the lift, drag, and pitching moment coefficients (as lumped parameters).

\subsection{Test case 2: autonomous underwater glider}
The second test case addresses the design optimization of a bio-inspired autonomous underwater glider (AUG) with a manta-like shape, that has been selected as a test case for the NATO-AVT-404 Research Task Group on ``Enhanced Design Processes of Military Vehicles through Machine Learning Methods''. The geometry is constructed as a continuous wing using a section-wise scheme and is by design spanwise symmetric; only half of the body is parameterized and then mirrored at the root section. As illustrated in Fig.~\ref{fig:glider_shapes}, the model comprises three zones: the center body, the transition region, and the outer wing. Four sections define the half-span: (i)~the root (section~1), (ii)~the end of the main body (section~2), (iii)~the end of the transition (section~3), and (iv)~the tip (section~4). Each section is fully determined by ten variables: four parameters for the section's NACA 4-digit airfoil (maximum camber \(m\), maximum camber position \(p\), thickness ratio \(t\), chord \(c\)); three parameters for positioning the section's leading edge \((x_0, z_0, s_0)\); three parameters defining the rotation angles (pitch $\vartheta$, roll $\phi$, yaw $\psi$), applied sequentially about the leading edge.
Section~1 retains only 2~degrees of freedom (\(t\) and \(c\)), while all other parameters are fixed. Consequently, the complete parametric model includes \(M=32\)~design variables. The manta-like geometry is constructed via the \texttt{OpenCASCADE} CAD kernel~\cite{OpenCascade} and meshed with \texttt{Gmsh} v4.12.1~\cite{geuzaine2009gmsh}.

The glider surface mesh consists of \(L=784\) elements, and the physical data (pressure coefficients, lift, and drag) are acquired via the \texttt{PUFFIn} solver~\cite{perali2024performance}, developed by ENSTA Bretagne, which combines an incompressible potential approach with viscous corrections. The operating conditions correspond to a uniform inflow velocity of \(0.25\,\text{m/s}\) in seawater at a depth of \(1500\,\text{m}\), and an 8-degree angle of attack that yields near-optimal efficiency for the baseline configuration.

\subsection{Test case 3: ship propeller}
\begin{figure}[!b]
    \centering
    \includegraphics[width=0.5\columnwidth]{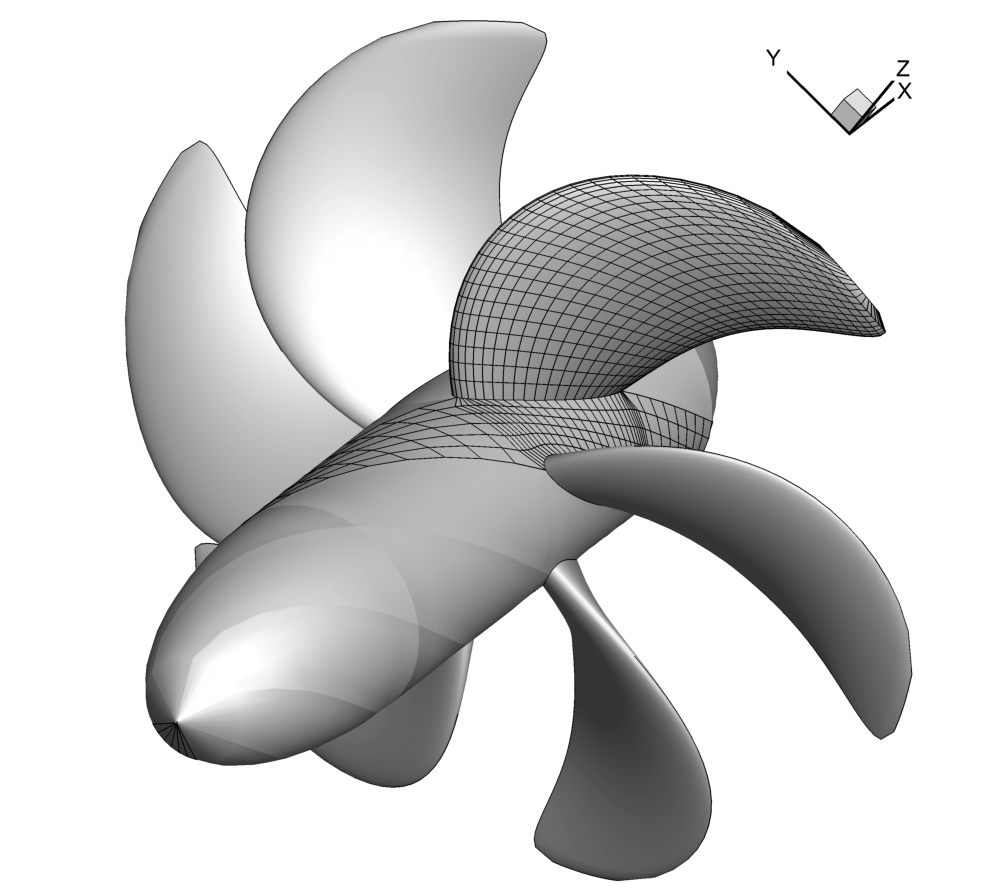}
   
    \caption{Test case~3: geometry of the reference ship propeller; ``Key-blade'' with the surface panel discretization used for BEM calculations of performances}
    \label{fig:3D_elica}
\end{figure}

\begin{figure}[!t]
    \centering
    \subfigure[Radial distributions (chord, pitch, rake, and maximum sectional camber, from left to right)]{
    \includegraphics[width=0.21\linewidth]{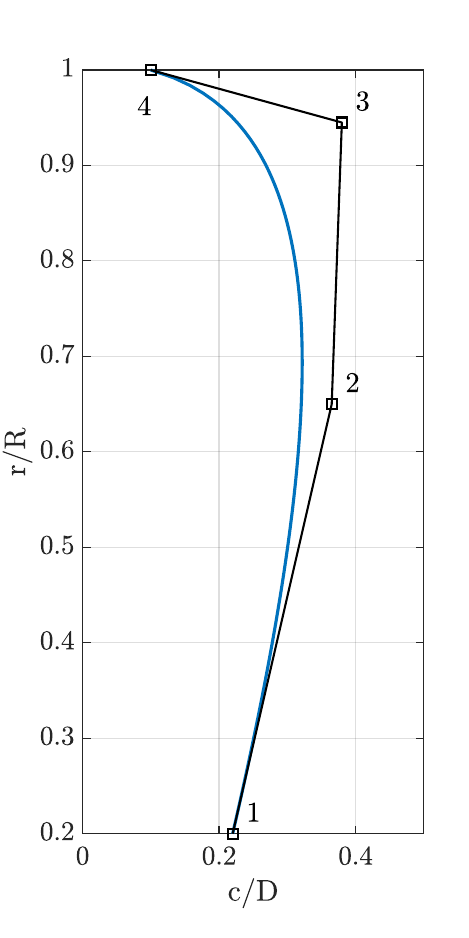}
    \includegraphics[width=0.21\linewidth]{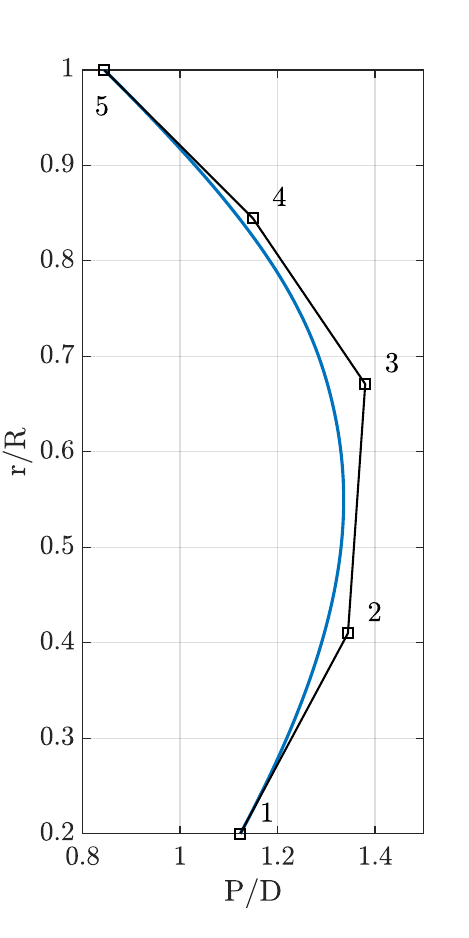}
    \includegraphics[width=0.21\linewidth]{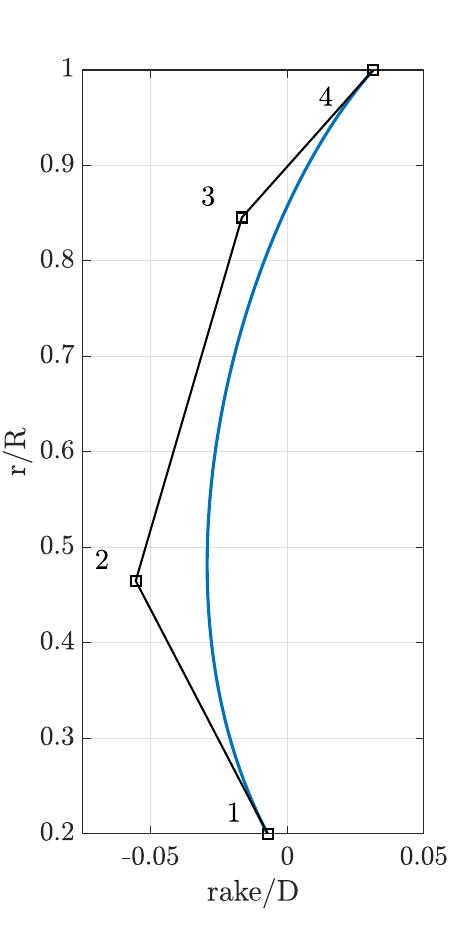}
    \includegraphics[width=0.21\linewidth]{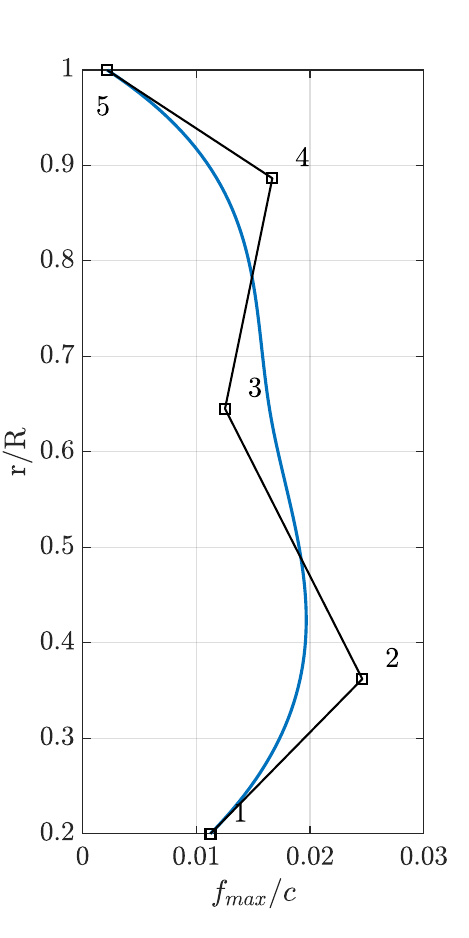}}    
\subfigure[Sectional distributions (chordwise thickness and camber line, left and right)]{
    \includegraphics[width=0.44\linewidth]{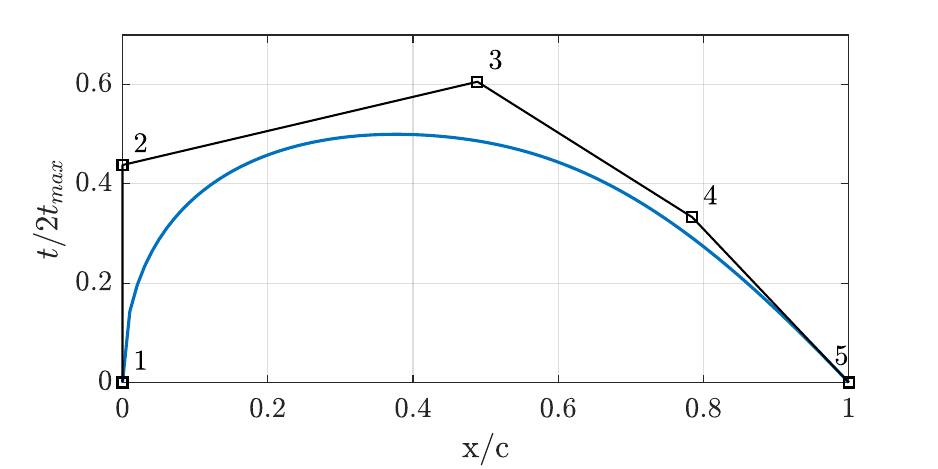}
    \includegraphics[width=0.44\linewidth]{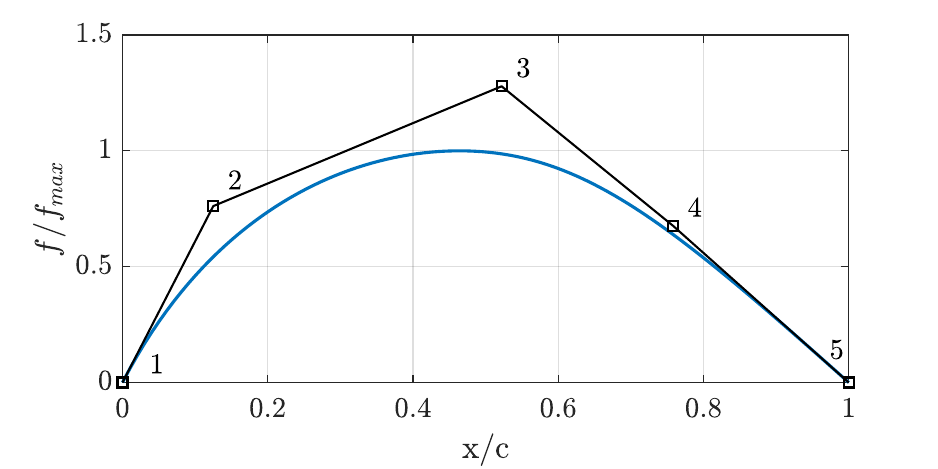}}
    \caption{Test case 3: distributions of the geometrical features of the propeller}
    \label{fig:prop_sezioni}
\end{figure}
Test case~3 considers the design optimization of a six-blade, right-handed marine propeller of a cruise ship. The reference propeller, shown in Fig.~\ref{fig:3D_elica}, was designed for a nominal advance coefficient ($J = V/nD$, where $V$ is the advance speed, $D$ the propeller diameter, and $n$ the rate of revolution in rps) of about 0.87 and a cavitation index $\sigma_n = 2(p-p_{vapor})/(\rho n^2D^2)$  of $2.25$. In this functioning condition, the delivered thrust corresponds to a thrust coefficient ($K_T = T/(\rho n^2D^4)$) of 0.186.
The parametric description consists of B-Spline curves that describe the geometrical features, in radial and chordwise directions, of the propeller blade. This approach, extensively validated in several design-by-optimization cases \cite{bertetta2012cpp, gaggero2020RIM}, provides a robust and easily controllable representation of the blade geometry in terms of its fairness. Moreover, it allows for the independent choice of quantities, like maximum sectional thickness, using criteria other than hydrodynamic shape optimization (in this case, structural strength of the blade) which instead are hardly controllable when using free-form deformations or B-Surfaces describing directly the suction and pressure side of the blade.

\begin{figure}[!t]
    \centering
    \includegraphics[width=0.5\columnwidth]{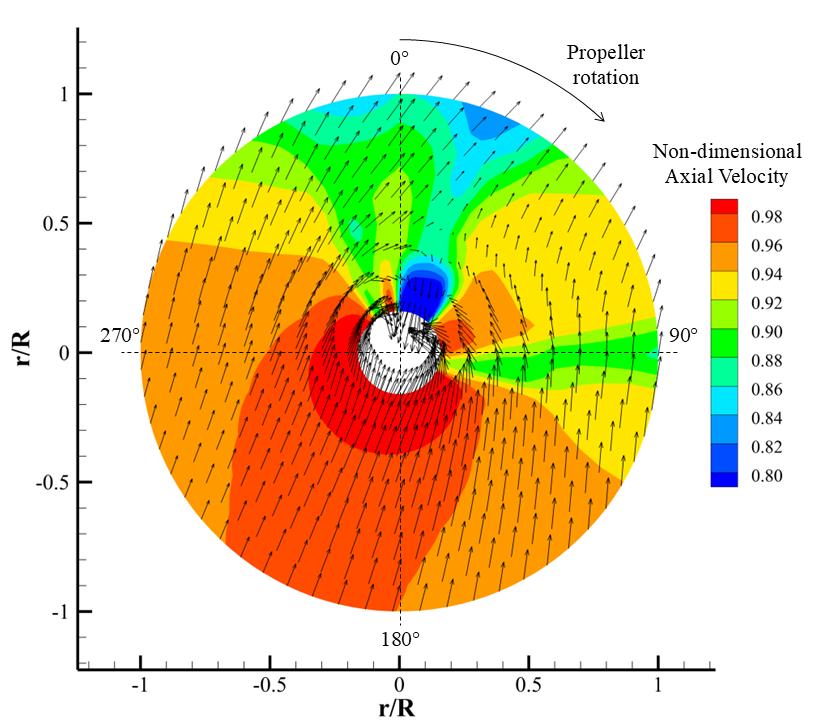}
   
    \caption{Non-uniform inflow (ship wake) to the propeller, seen from aft}
    \label{fig:Ship_Wake}
\end{figure}

For this particular problem, the blade geometry is parametrized through B-Spline control polygons defining radial distributions of non-dimensional chord ($c/D$, parameter name $c$), pitch ($P/D$, parameter name $pd$), rake ($rake/D$, parameter name $r$), and maximum sectional camber ($f_{max}/c$, parameter name $f$) of Fig.~\ref{fig:prop_sezioni}a. The shape of sectional hydrofoils is included as well in the parametric model by describing, again by using appropriate control polygons (see Fig.~\ref{fig:prop_sezioni}b), the non-dimensional chordwise sectional thickness ($t/2t_{max}$, parameter name $th$) and camber ($f/f_{max}$, parameter name $fh$) distributions of the blade profile which, radial section per radial section, are finally scaled with the corresponding maximum values. Given some constraints (the radial position of control points at root and tip is fixed, leading- and trailing-edge points are always given, and the chord at the tip never changes to comply with some limits of the flow solver adopted for the hydrodynamic characterization of performances), control points, sequentially numbered from root to tip (or from leading- to trailing-edge), are free to move within assigned ranges. Subscript ``$_x$'' indicates the radial, or chordwise, modification of the control point position, and subscript ``$_y$'' refers instead to changes in the quantity the control point describes. This parametrization leads to a total of \(M = 38\) design variables (5 for the chord, 8 for the maximum sectional camber and pitch each, 6 for the rake, 6 for the sectional camber line, and 5 for the sectional thickness non-dimensional distribution) which describe entirely the blade shape with the exception of the skew that is maintained unchanged and identical to that of the reference propeller.   
For the hydrodynamic analyses, a low-fidelity boundary element method (BEM) is employed. The BEM code  \cite{Gaggero_BEM} was developed at the University of Genoa since early 2000 for analysis and design-by-optimization purposes. It is a Morino, Dirichlet-type boundary condition implementation of panel methods for incompressible, potential flow solution, which makes use of the ``key blade'' approach to deal with stationary or unsteady problems \cite{kinnas1992boundary}. It includes a cavitation model (the sheet cavity model at leading edge and midchord, both on suction and pressure side proposed in \cite{fine1992nonlinear}), wake alignment capabilities and the iterative Kutta condition.
For current analyses, the ``key blade'' is discretized with 1250 hyperboloidal panels (\L = 1326\ nodes). The collection of physical information was carried out by considering several equivalent steady-state operating conditions representative of the most critical conditions encountered by the blade during a revolution, given the design of the reference propeller, which was developed to operate in the spatially non-uniform wake of the ship shown in Fig.~\ref{fig:Ship_Wake}. 

\begin{figure*}[!b]
    \centering
    \includegraphics[width=0.8\linewidth]{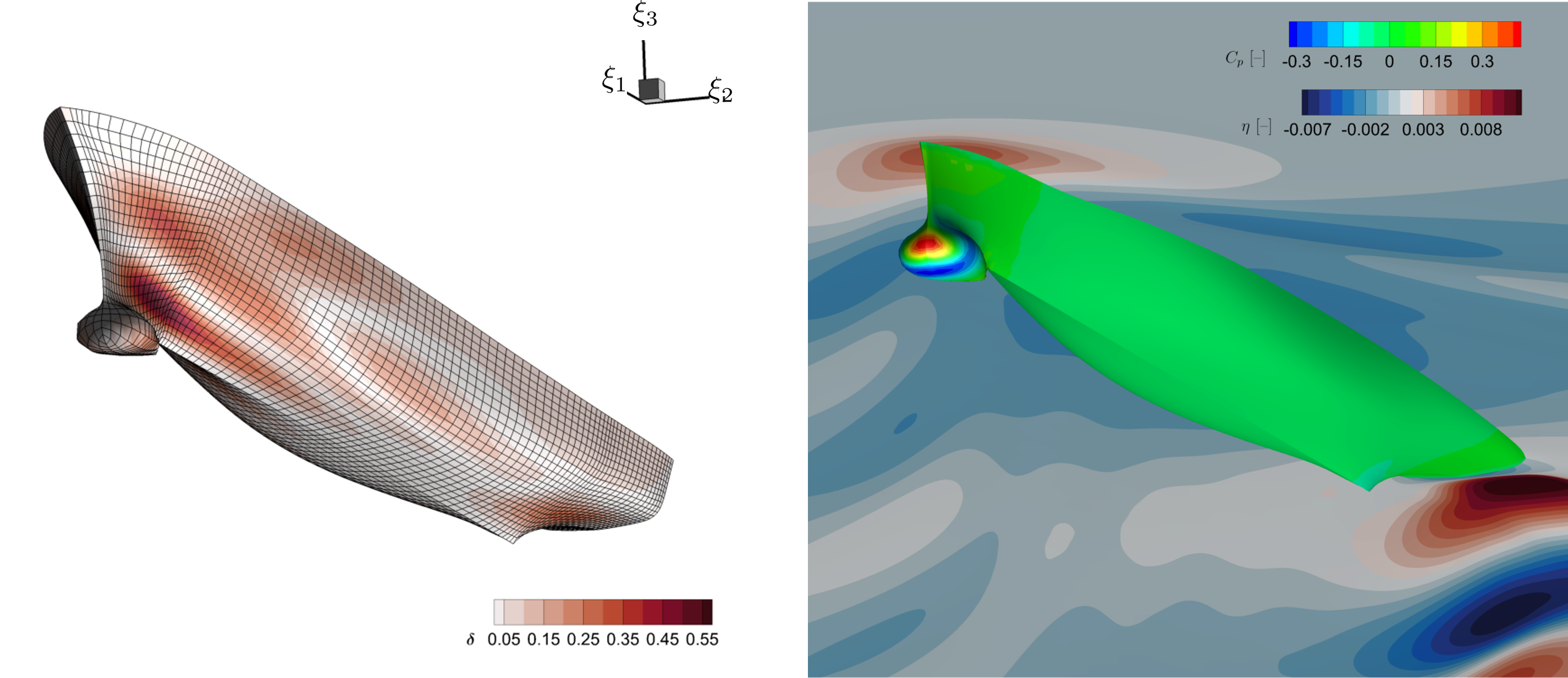}
    \caption{Test case~4: Example of (left) shape modification magnitude with discretization and (right) calm-water numerical solution for the pressure distribution and the wave elevation pattern $\eta$ generated by the 5415 hull}
    \label{fig:5415_shapes}
\end{figure*}

Despite the computational efficiency of the Boundary Element method, addressing the truly unsteady performance of the propeller is excessively demanding for the characterization of the thousands of different configurations needed to feed the physically-informed parametric model embedding methods. Conversely, steady analyses under equivalent flow conditions have consistently yielded the performance indicators of conventional design tools \cite{lerbs1952moderately, epps2013unified} and have been effectively utilized in numerous designs by optimization processes \cite{gaggero2025robust} as a surrogate for unsteady functioning predictions. In addition to the equivalent nominal functioning condition (i.e., the propeller operating the circumferentially averaged axial inflow wake), two additional calculations were included in the analyses to provide a broader view of the propeller performance. One corresponds to the loaded case of the blade passing through the 90° position of the non-uniform inflow wake, where the action of tangential velocities increases the angle of attack. The other, on the contrary, addresses the unloaded blade when at 270° into the wake, where instead the tangential velocities act to reduce the angle of attack. Together, they allow collecting distributed pressure data as well as lumped performance indicators such as thrust coefficient $K_T$, efficiency $\eta$, and tip-vortex intensity $\Gamma$ related to the risk of different types of cavitation (suction or pressure side, respectively loaded and unloaded equivalent conditions) and to expected unsteady behaviour of the propeller.

\subsection{Test case 4: 5415 destroyer-type vessel}
The fourth test case targets the international benchmark 5415 hull \cite{serani2024hydrodynamic}, which is a geosymmetric replica of the DDG-51 (a US Navy destroyer-type vessel). Shape parameterization is based on a recursive set of \(M=27\) global modification functions applied within a hyper-rectangular region enclosing the demi-hull~\cite{serani2016ship,serani2022hull}. Figure \ref{fig:5415_shapes} (left) shows an example of the shape modification magnitude for one design variant.

The surface mesh of the hull comprises \(L=2250\) grid nodes. Physical data include calm-water and seakeeping performance, extracted from the in-house linear potential-flow solver \texttt{WARP}~\cite{bassanini1994-SMI} (developed by CNR-INM) with viscous corrections, alongside a linearized strip-theory code \texttt{SMP}~\cite{meyers1985-TechRep}. Calm-water analyses provide the distribution of surface pressure (see Fig. \ref{fig:5415_shapes}, right) and the wave-resistance coefficient; seakeeping analyses deliver estimates of the pitch motion $(\theta)$ RMS and the vertical acceleration $(a_B)$ RMS at the bridge. These integrated and distributed quantities collectively serve as physical inputs for PI-PME and PD-PME, allowing for a thorough exploration of how variations in hull geometry correlate with performance in real sea conditions.

%%%%%%%%%%%%%%%%
\begin{figure}[!b]
    \centering
    \subfigure[TC\#1: RAE-2822 ]{\includegraphics[width=0.49\linewidth]{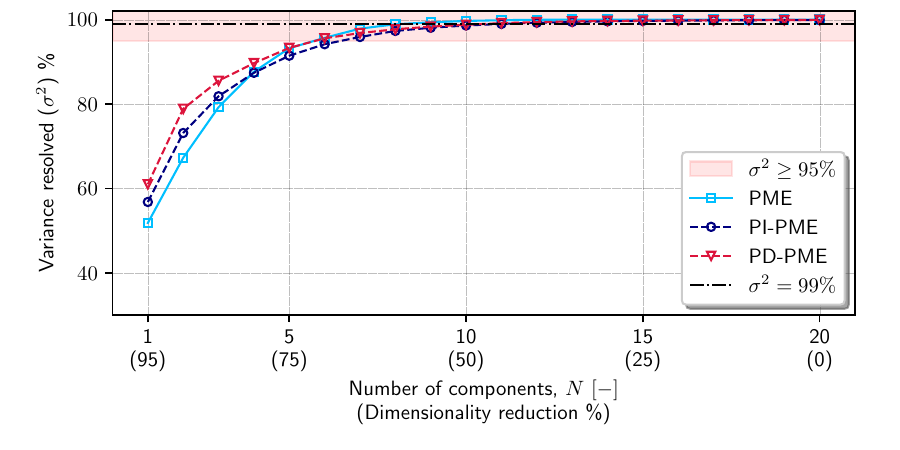}}
    \subfigure[TC\#2: AUG      ]{\includegraphics[width=0.49\linewidth]{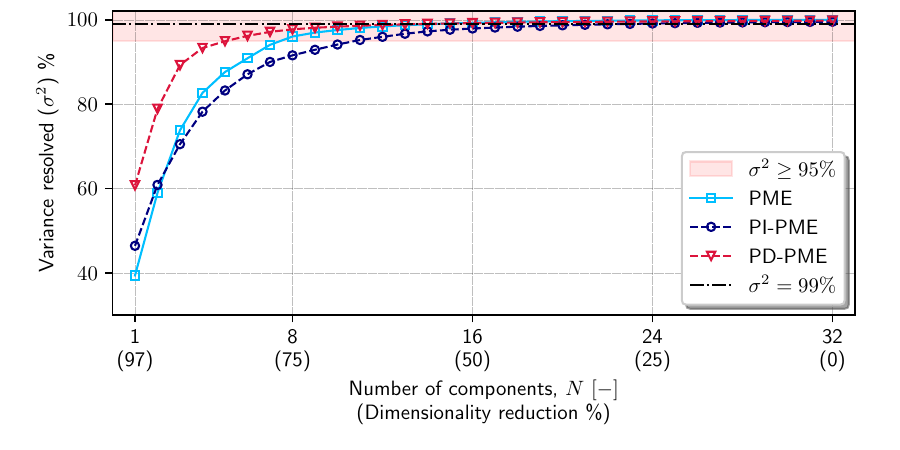}}\\
    \subfigure[TC\#3: Propeller]{\includegraphics[width=0.49\linewidth]{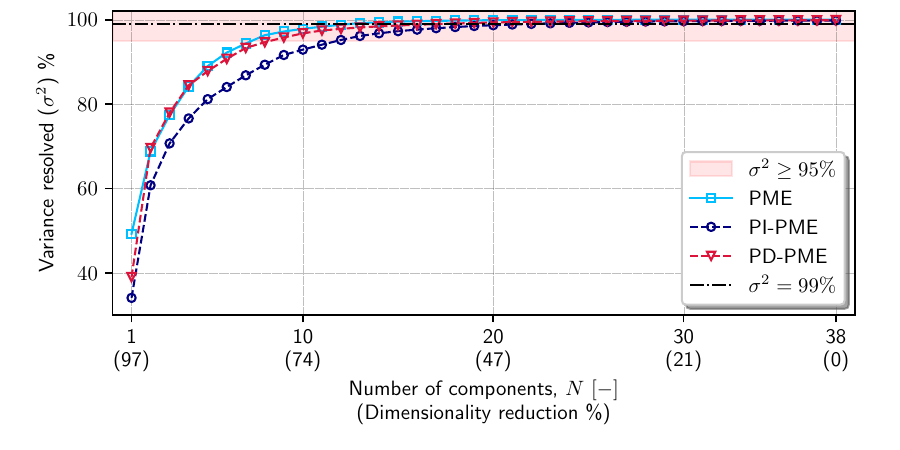}}
    \subfigure[TC\#4: 5415     ]{\includegraphics[width=0.49\linewidth]{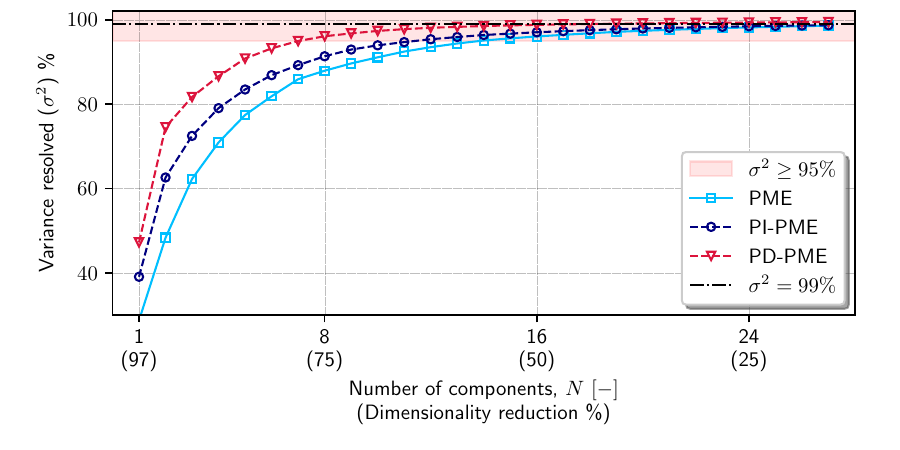}}
    \caption{Variance resolved as a function of the reduced design variables}
    \label{fig:eigsums}
\end{figure}
\section{Results}
To train the standard PME and its physics-informed extensions (PI-PME and PD-PME), an initial set of $S=16385$ samples was generated for each test case. \change{To ensure a uniform and representative coverage of the design space, the initial dataset is generated using a Sobol low-discrepancy sequence. This approach provides better space-filling properties than purely random sampling and avoids clustering, which is critical to reliably constructing the reduced manifold.}
In addition to covering the geometric design space, physical simulations were performed at these samples whenever possible, providing distributed or lumped physical data. Any geometries that produced unfeasible outcomes (e.g., degenerate shapes, solver divergence, or \texttt{NaN} outputs) were discarded. This filtering process alone already serves as a preliminary assessment of the chosen parameterization, as it highlights whether certain geometric variations are physically inadmissible or numerically unstable.

An additional layer of screening was then applied to exclude extreme outliers. Specifically, all geometric--physical samples were retained only if the corresponding physical outputs fell within the ranges
$[Q1-3IQR,Q3+3IQR]$, where $Q1$ and $Q3$ are the first and third quartiles, respectively, and $IQR=Q3-Q1$. 
\change{This distribution-agnostic rule was selected to balance permissiveness and robustness: it allows broad variability across designs while excluding pathological outliers.} As a result, only physically credible shapes and simulations populate the final data set used for dimensionality reduction.
\change{The filtering strategy is intentionally cautious, aiming to eliminate anomalous or ill-posed configurations that could contaminate the reduced basis with non-representative trends. It is finally worth noting that, while the reduced space may still contain directions that could lead to poorly behaved geometries, standard constraint-handling techniques in the subsequent optimization process can be employed to reject such configurations if encountered.}

\begin{table}[!t]
    \centering
        \caption{Design-space dimensionality reduction summary for retaining 95\% of the problem variance}
    \begin{tabular}{ccccccc}
    \toprule
    TC & Geometry & $S$ & $M$         &  \multicolumn{3}{c}{$N$ (dimensionality reduction\%)}\\
     \#& model    & Samples & Original    &  PME & PI-PME & PD-PME\\
         \midrule
    1    & RAE-2822  & 6,022 & 20  &  6 (70\%)&  7 (65\%) &  6 (70\%)\\
    2    & AUG       & 5,330 & 32  &  8 (75\%)& 11 (66\%) &  5 (84\%)\\
    3    & Propeller &13,615 & 38  &  8 (79\%)& 12 (68\%) &  8 (79\%)\\
    4    & 5415      & 5,842 & 27  & 14 (48\%)& 12 (66\%) &  7 (74\%)\\
         \bottomrule
    \end{tabular}
    \label{tab:dr_summary}
\end{table}

\begin{figure*}[!t]
    \centering
    \includegraphics[width=0.49\linewidth]{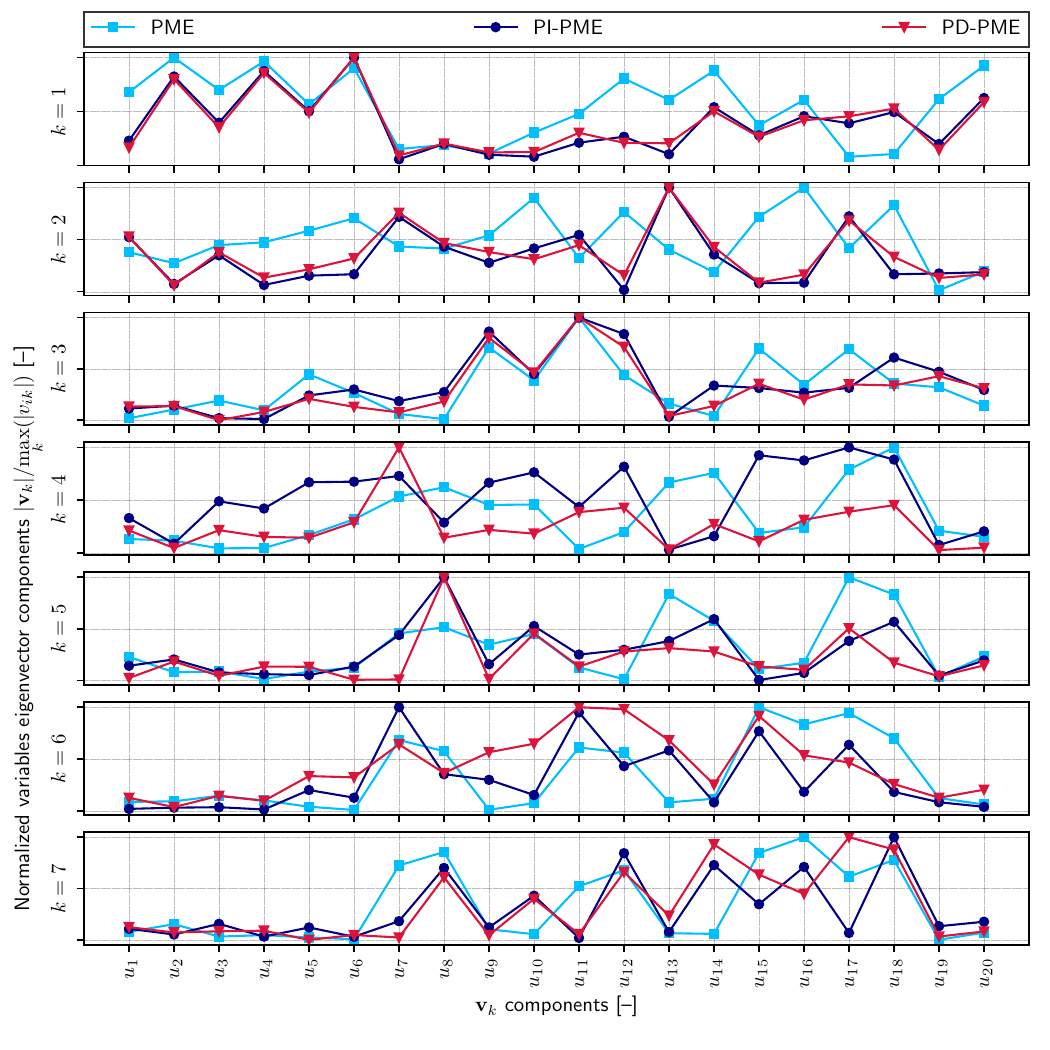}
    \includegraphics[width=0.49\linewidth]{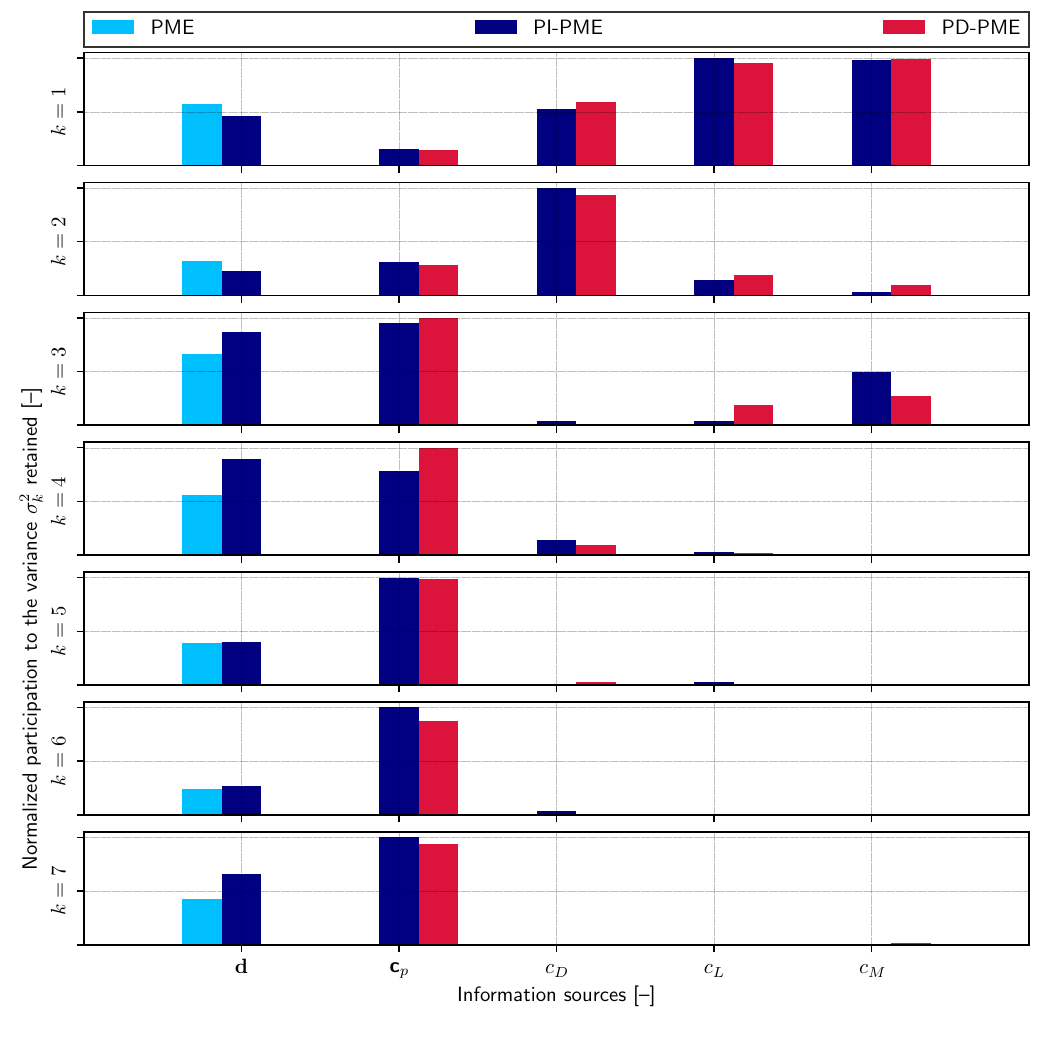}
    \caption{Design-space dimensionality reduction (left) eigenvectors $\bfv_k$ that embed the original design variables and (right) participation to the variance retained by each eigenvector of geometrical and physical information for RAE-2822 (TC\#1)}
    \label{fig:rae_eigs}
\end{figure*}

\begin{figure*}[!t]
    \centering
    \includegraphics[width=0.49\linewidth]{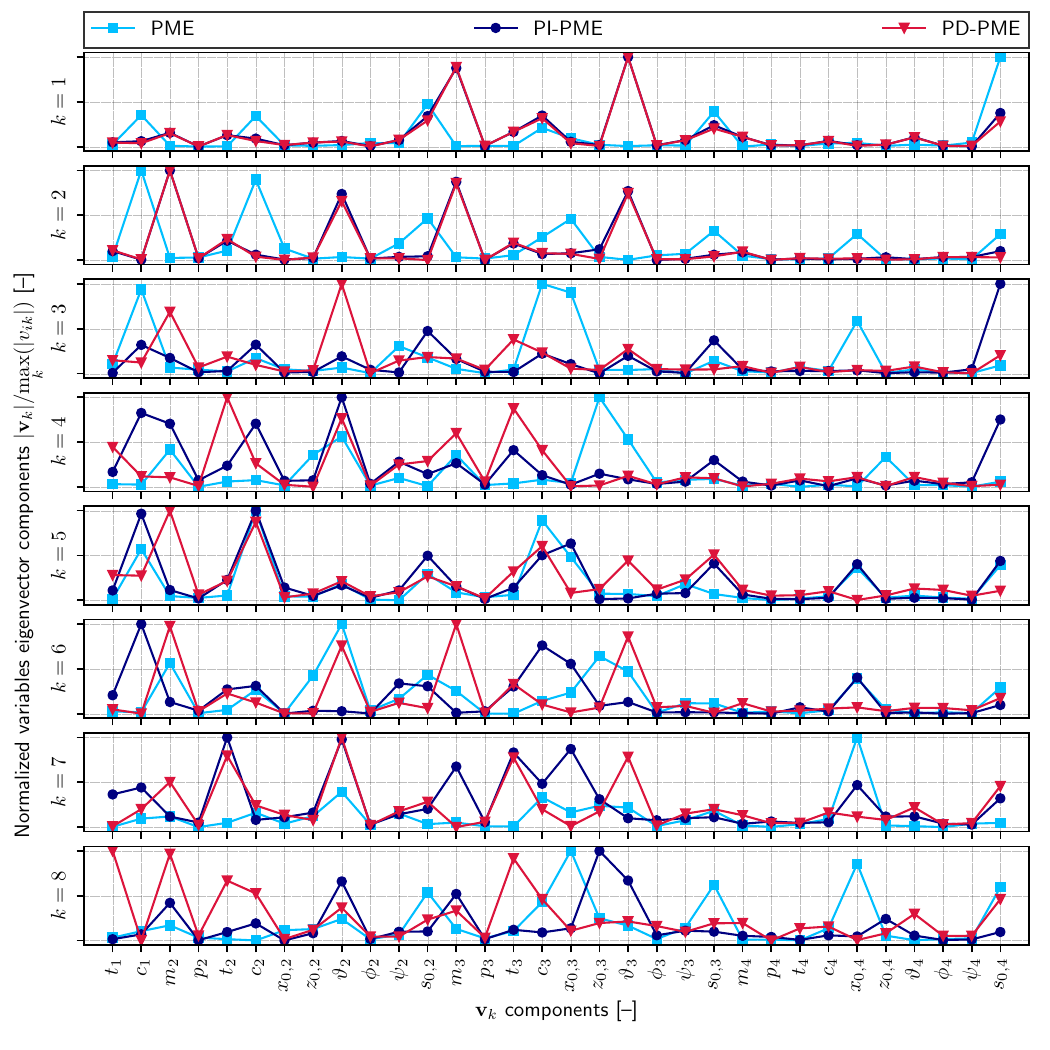}
    \includegraphics[width=0.49\linewidth]{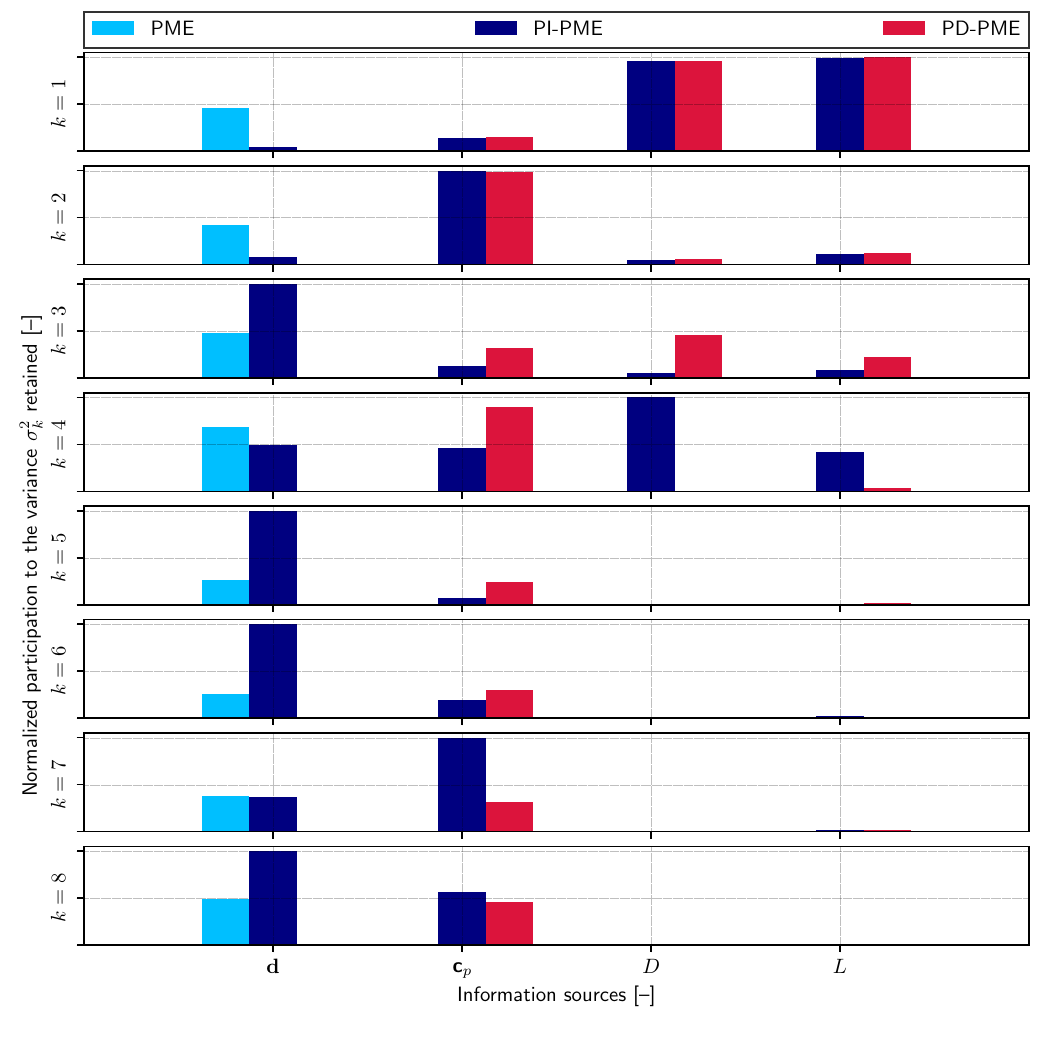}
    \caption{Design-space dimensionality reduction (left) eigenvectors $\bfv_k$ that embed the original design variables and (right) participation to the variance retained by each eigenvector of geometrical and physical information for AUG (TC\#2)}
    \label{fig:glider_eigs}
\end{figure*}

\begin{figure*}[!t]
    \centering
    \includegraphics[width=0.49\linewidth]{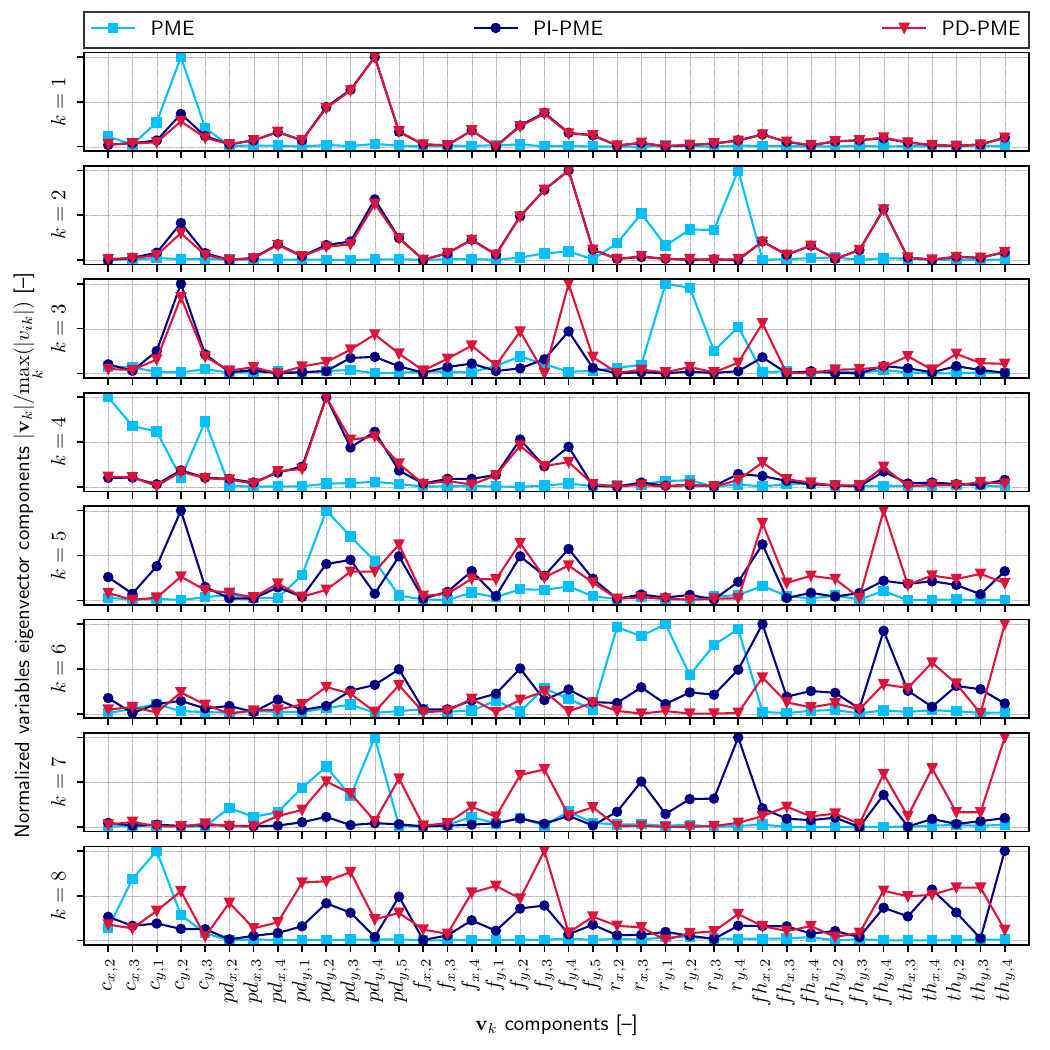}
    \includegraphics[width=0.49\linewidth]{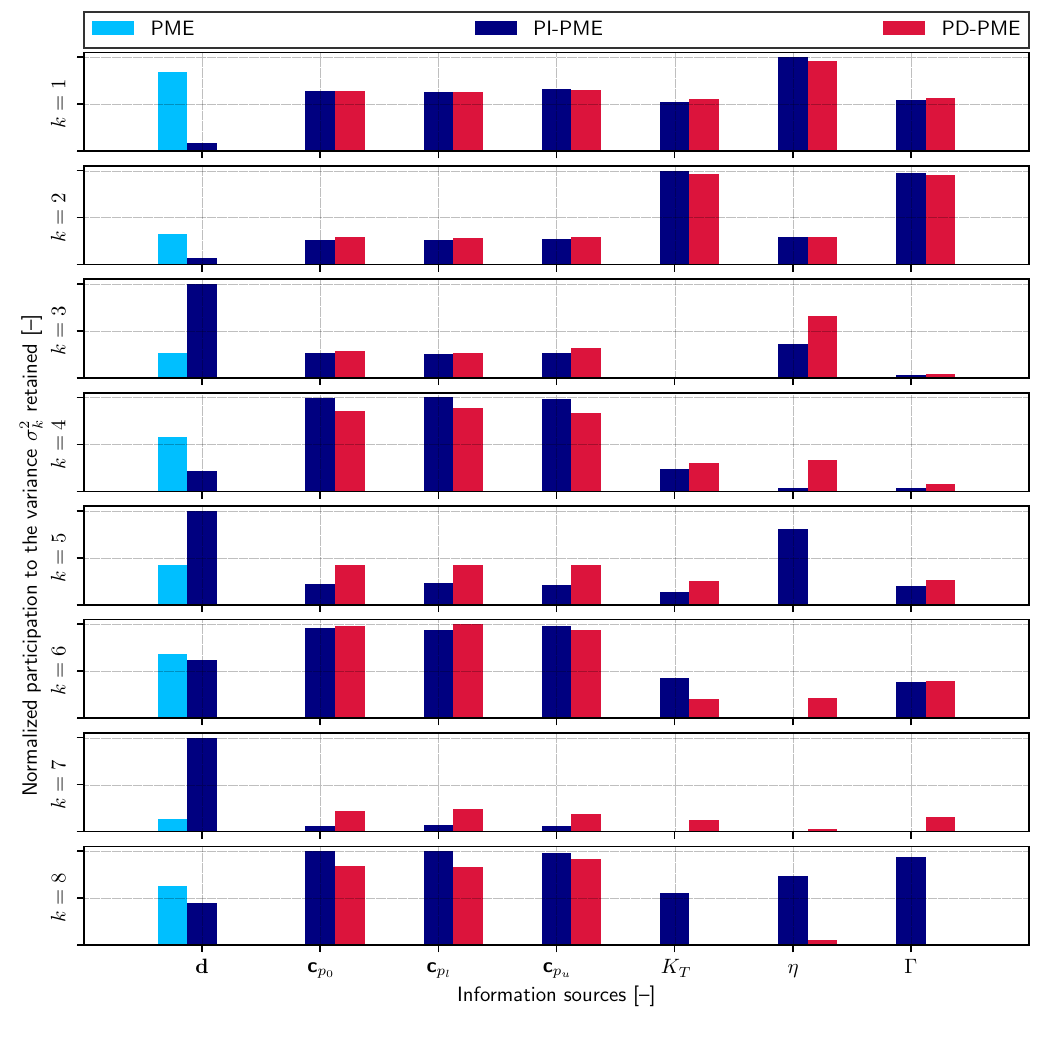}
    \caption{Design-space dimensionality reduction (left) eigenvectors $\bfv_k$ that embed the original design variables and (right) participation to the variance retained by each eigenvector of geometrical and physical information for propeller (TC\#3)}
    \label{fig:propeller_eigs}
\end{figure*}

\begin{figure*}[!t]
    \centering
    \includegraphics[width=0.49\linewidth]{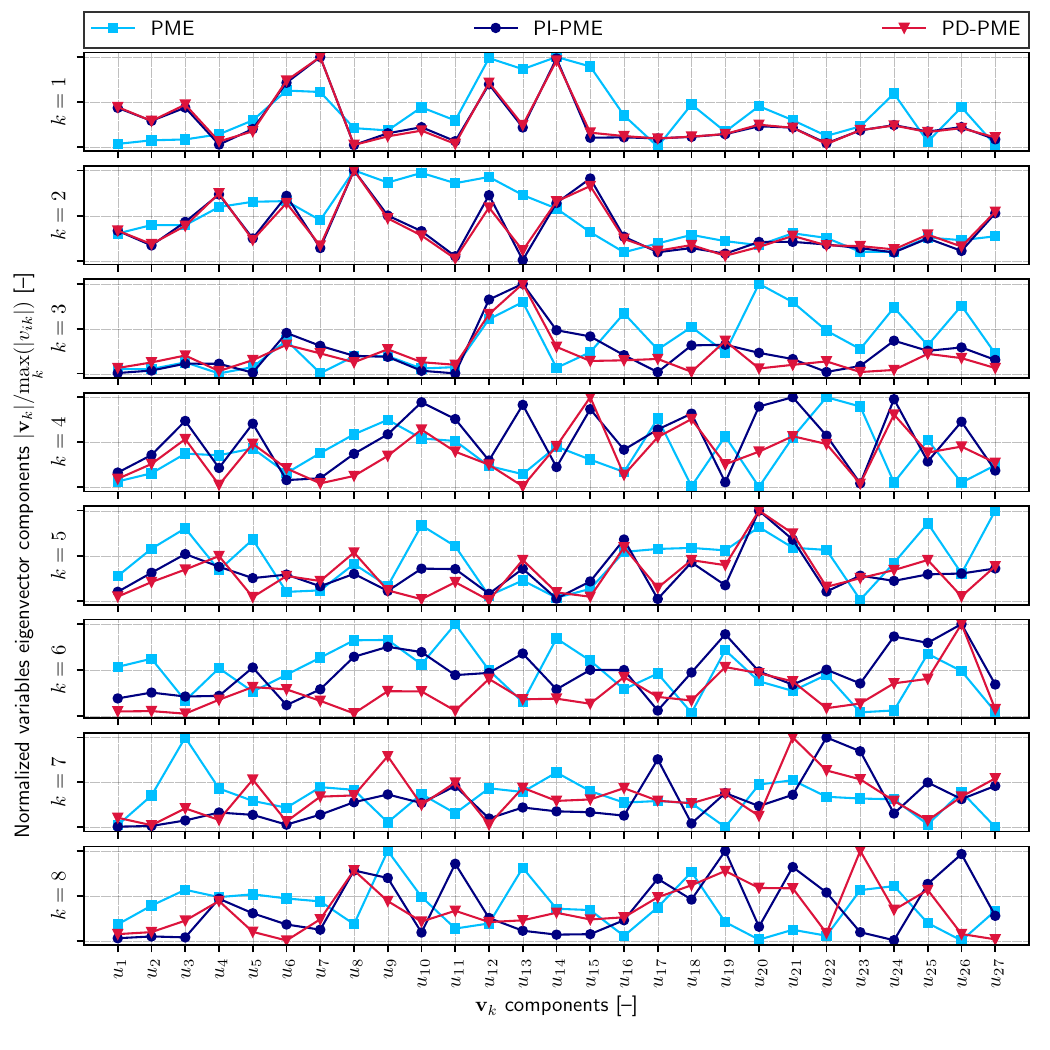}
    \includegraphics[width=0.49\linewidth]{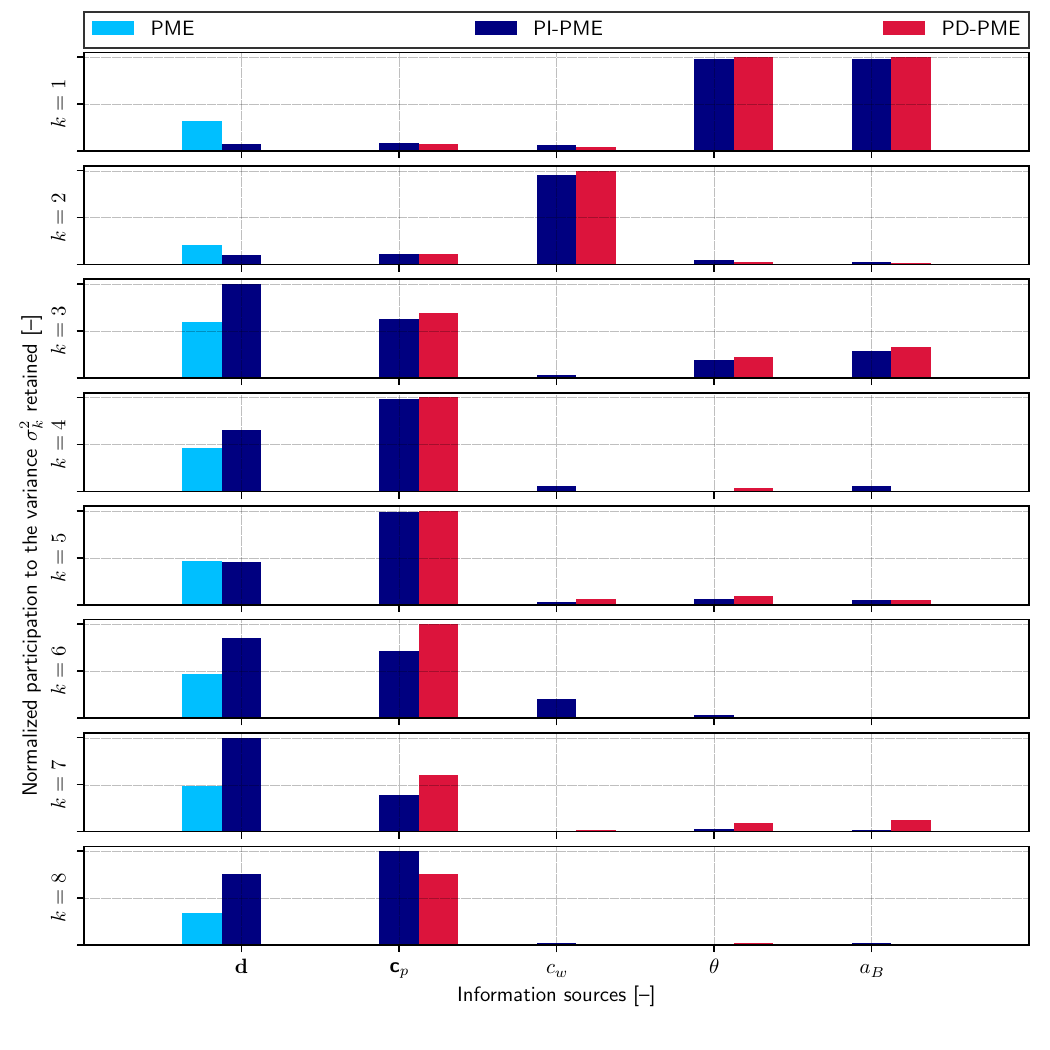}
    \caption{Design-space dimensionality reduction (left) eigenvectors $\bfv_k$ that embed the original design variables and (right) participation to the variance retained by each eigenvector of geometrical and physical information for 5415 (TC\#4)}
    \label{fig:5415_eigs}
\end{figure*}

Figure \ref{fig:eigsums} presents how the fraction of variance ($\sigma^2$) evolves as a function of the number of reduced components $N$, contrasting PME (purely geometric variance), PI-PME (geometry + physical data), and PD-PME (physical data only). These results are numerically summarized in Table \ref{tab:dr_summary} (e.g., how many modes are needed to reach 95\% of the total variance). Meanwhile, Figures \ref{fig:rae_eigs}–\ref{fig:5415_eigs} provide deeper insights into how that variance is organized for each mode: on the left, each figure shows the normalized eigenvector components $\vert\mathbf{v}_k \vert /\underset{k}{\max}(\vert v_{ik}\vert)$ along the vertical axis, plotted across the original design variables on the horizontal axis. Large peaks indicate which design variables dominate that particular mode; on the right is shown the participation of each data source (e.g., geometry $\mathbf{d}$, distributed pressure $\mathbf{c}_p$, or lumped coefficient like $c_D$, $c_L$, etc.) to variance retained by each mode. The latter plots reveal whether a mode primarily represents geometric variation, physical variation, or a mixture of both.

Taken together, these figures clarify whether geometry and physics are strongly correlated (leading to fewer modes or shared modes) or largely uncorrelated (leading to more modes, each distinctly owned by geometry or physics). Below is a case-by-case discussion.

PME, PI-PME, and PD-PME curves in the RAE airfoil problem lie close together (see Fig. \ref{fig:eigsums}a); each method needs a comparable number of modes ($N=6$ or 7) to reach 95\% or 99\% variance. This suggests that geometric changes coincide strongly with variations in the pressure distribution and integrated aerodynamic coefficients.
Figure \ref{fig:rae_eigs} (left) shows a pretty similar embedding with eigenvectors almost coincident for PME, PI-PME, and PD-PME (see first and third mode), in particular, it can be seen how the first mode is mainly participated by $u_2$, $u_4$, and $u_6$, that correspond to the polynomial shape modification function applied to the lower side of the airfoil, highlighting also how the lower side modification is the major region affecting the aerodynamic performance. This is further confirmed by Figure \ref{fig:rae_eigs} (right) where is shown that geometry ($\bfd$) and physical parameters ($\mathbf{c}_p$, $c_L$, $c_D$, and $c_M$) share most of the first modes: apart from mode 2, that is mainly participated by the drag coefficient, no distinct pure physics or pure geometry mode dominates the decomposition.
Because airfoil shape changes strongly affect the aerodynamic response under these conditions, the geometry and physics are highly correlated. Hence, adding physical data (PI-PME) does not inflate the necessary dimensionality, and removing geometry altogether (PD-PME) does not drastically reduce the dimension beyond PME, meaning that both geometry and physics point in similar directions in the design space.

For the glider test case, the PME rapidly achieves a high fraction of geometric variance ($N=8$ modes for 95\%), whereas the PI-PME requires more modes ($N=11$). In contrast, PD-PME dips below the PME line, reaching 95\% in fewer modes ($N=5$), as shown in Fig. \ref{fig:eigsums}b. Figure \ref{fig:glider_eigs}, shows, on the left, that while for PME the highest variance (retained by the first mode) is mainly due to the span ($s_{0,4}$) of the tip section, the variance of PI-PME and PD-PME are mainly affected by the maximum camber and the pitch angle of the end of the transition from the center body to the outer wing ($m_3$ and $\vartheta_3$). As shown in Fig. \ref{fig:glider_eigs} (right) the first modes appear ``decoupled'' (i.e., geometry heavy or physics heavy): for both PI-PME and PD-PME the first mode is mainly participated by global forces ($L$  and $D$), the second mainly by the pressure distribution $\mathbf{c}_p$, and the third mainly by the geometry (for PI-PME only).
There is a partial uncorrelation between geometry and physics: large shape modifications do not necessarily produce proportional changes in lift or drag. Hence, PI-PME must retain more overall modes to capture both geometry-dominated and physics-dominated directions. Meanwhile, if the goal is purely performance-driven, PD-PME focuses on those physical variations alone, achieving a stronger reduction with fewer modes.

Concerning the propeller test case, PME captures geometric variability well in about $N=8$ modes, but PI-PME demands significantly more ($N=12$) to cover combined variance from geometry plus thrust, efficiency, vortex data, and pressure distributions. PD-PME, by contrast, can exceed 95\% in only $N=8$ modes if focusing purely on physical data. Figure \ref{fig:propeller_eigs} shows similar results as for the glider. 
% Stefano
The first mode of PME is participated by the variation of the intermediate control point of the chord ($c_{y,2}$), which obviously causes the most relevant modifications of the propeller shape, interpretable as variations of the expanded area of the blades. The second and third modes of the PME identify the rake (through $r_{y,2}$ and $r_{y,4}$) as the second most influential parameter on the geometry variance. This behavior is also plausible, as the ``shift'' of the blade sections in the longitudinal plane, which is the modification induced by the rake, contributes more to the differences in the shapes of the blade than the localized modifications induced by changes in camber and thickness (by radial or sectional perturbations of control points).  
%Stefano
The first two modes of PI-PME and PD-PME are mainly participated by the physical information, while the third of PI-PME is mainly participated by the geometry, even if in this case a partial participation of the physical quantities is visible.
% Stefano
In this respect, it is interesting to observe the very nice correlations between geometrical features and physical information revealed by both PI-PME and PD-PME. The first mode for both approaches, which is mainly participated by efficiency, identifies the pitch, and in particular its values at the tip, as the most important factor responsible for the mode, which is a trend exactly in line with traditional design method outputs and guidelines. The second mode, for which again PI-PME and PD-PME are perfectly overlapped, establishes a correlation between the delivered thrust and the strength of the tip vortex circulation with the pitch and maximum values of camber again, which are, by analogy with 2D hydrofoils, responsible for sectional forces (lift) through the angle of attack (pitch) and lift coefficient at zero angle of attack (camber). When, instead, the focus is on physical information like the pressure coefficient, physics-informed models highlight the influence of pitch (responsible for hydrofoil suction peaks at the leading edge) and camber (responsible for pressure distributions on the central part of the hydrofoil) through modes 4 to 8. When purely geometrical analyses were employed, quantities like maximum sectional camber and hydrofoil sectional shape (camber and thickness) were completely discarded, being, from a geometrical point of view, responsible for minimal and localized shape variations, in favor of the rake that is instead responsible for large geometrical modification but has almost no influence on (any) propeller performance indicators. Only embedding the physical quantities activates these features, which, indeed, regardless of the small geometrical variations they induce, have a critical influence on pressure distributions and the overall performance of the propeller.
% Stefano
Geometry and physical performance are weakly correlated. This necessitates many extra modes in PI-PME to accommodate both sets of unaligned directions. PD-PME proves highly compressible because the pure physics space (i.e., lumps + distributions) is lower-dimensional in meaningful directions for performance.

Finally, looking at the 5415 results (see Fig. \ref{fig:5415_eigs}, PME needs $N=14$ modes to exceed 95\% variance, PI-PME, interestingly, needs fewer ($N=12$) modes to reach the same coverage, while PD-PME reaches 95\% variance in only $N=7$ modes. The first mode in PI-PME is driven almost exclusively by lumped seakeeping parameters. The second mode largely captures the calm-water wave-resistance coefficient. Starting from the third mode onward, geometry and distributed pressure tend to move ``in lockstep,'' indicating a tight correlation between shape changes and local pressure distribution. 
At first glance, one might expect adding physics increases the total variance, as per the other three test cases. However, the crucial factor is that many purely geometric variations do not correlate with changes in both calm-water and seakeeping physical quantities. When geometry is combined with physics in PI-PME, the modes that represent ``empty'' geometric variance (irrelevant to performance) no longer contribute significantly to the principal components. As a result, fewer total modes can end up describing all relevant variations.

\begin{figure*}[!t]
    \centering
    \subfigure[TC\#1: RAE-2822 ]{\includegraphics[width=0.245\linewidth]{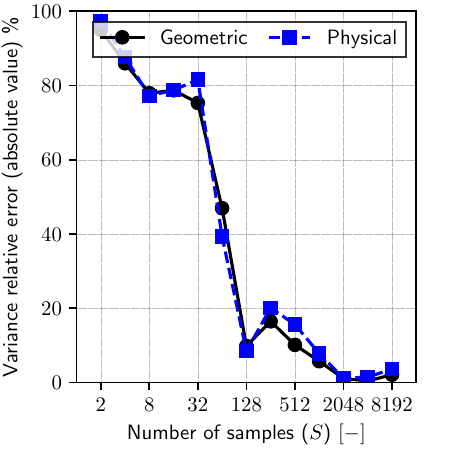}}
    \subfigure[TC\#2: AUG ]{\includegraphics[width=0.245\linewidth]{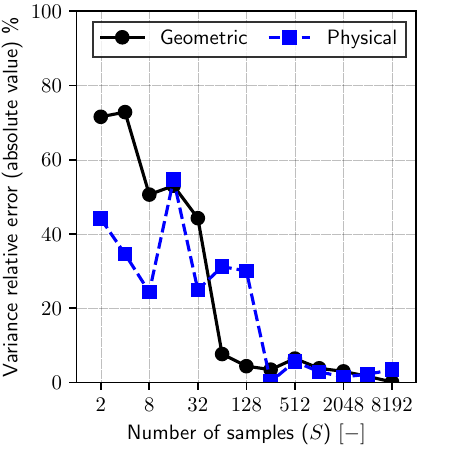}}
    \subfigure[TC\#3: Propeller ]{\includegraphics[width=0.245\linewidth]{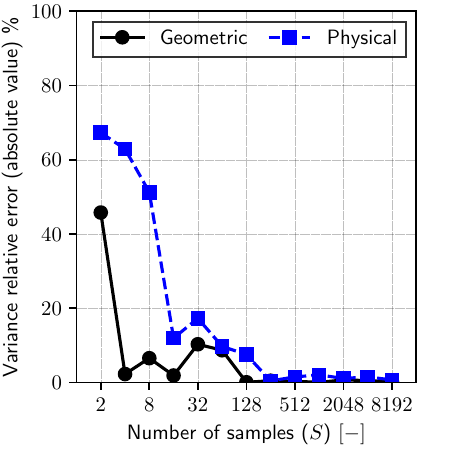}}
    \subfigure[TC\#4: 5415 ]{\includegraphics[width=0.245\linewidth]{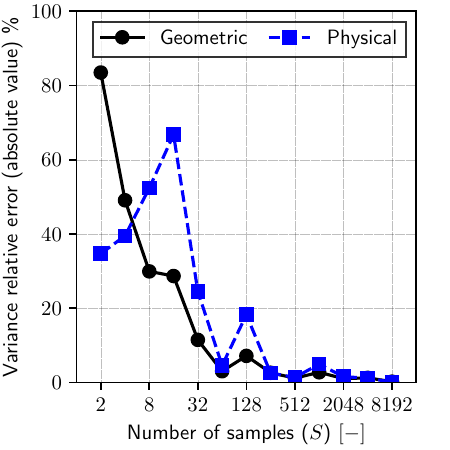}}
    \caption{\change{Geometric and physical variance sensitivity to the number of samples. The curves show the relative error compared to the maximum sample size (16,384), used as reference.}}
    \label{fig:var_conv}
\end{figure*}
%

%\newpage
\begin{figure}[!b]
    \centering
    \subfigure[TC\#1: RAE-2822]{
    \includegraphics[width=0.31\linewidth]{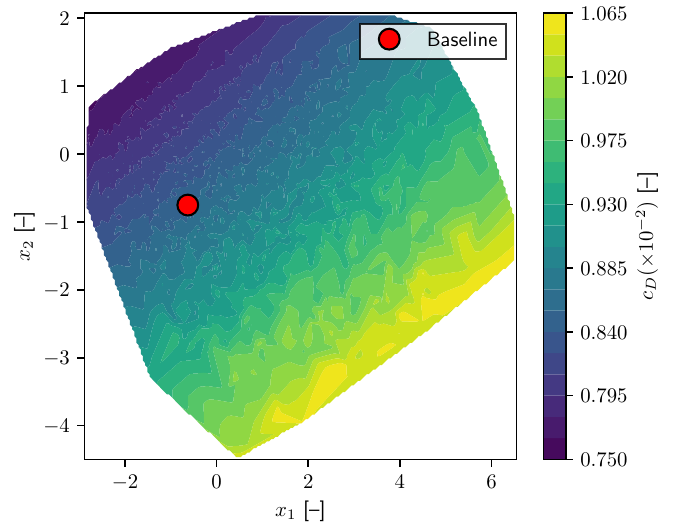}
    \includegraphics[width=0.31\linewidth]{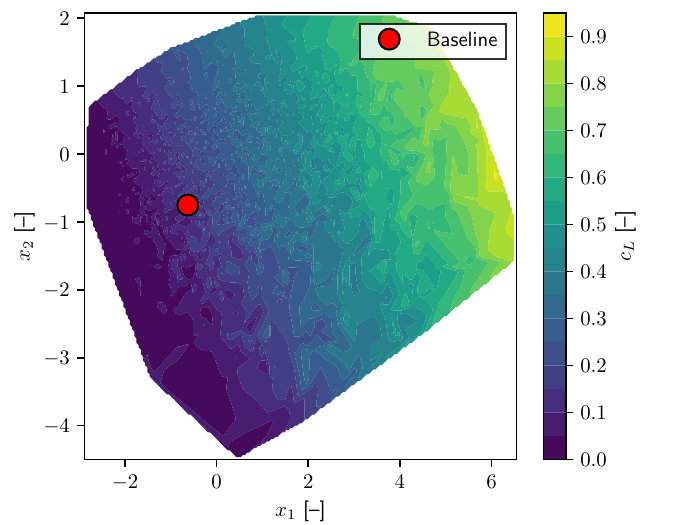}
    \includegraphics[width=0.31\linewidth]{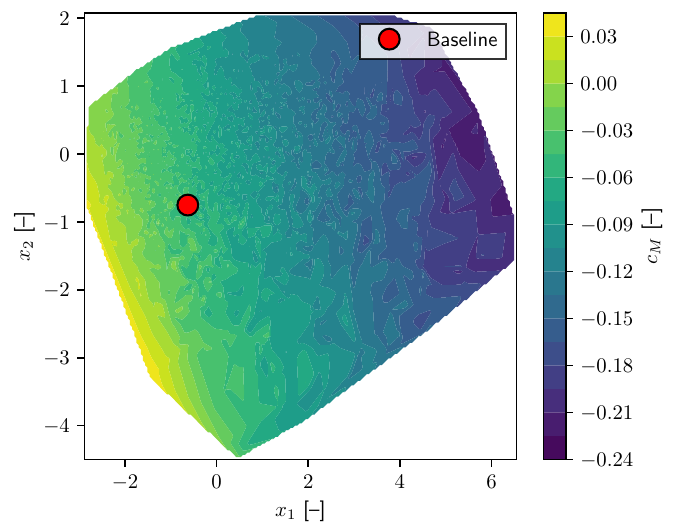}}
    \subfigure[TC\#2: AUG]{
    \includegraphics[width=0.31\linewidth]{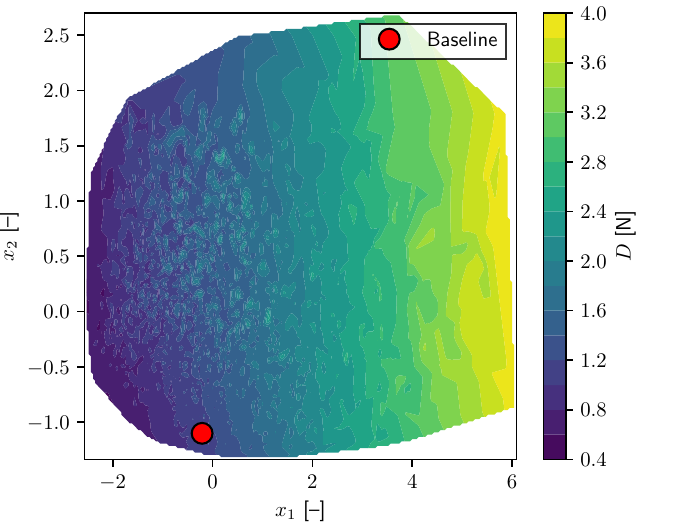}
    \includegraphics[width=0.31\linewidth]{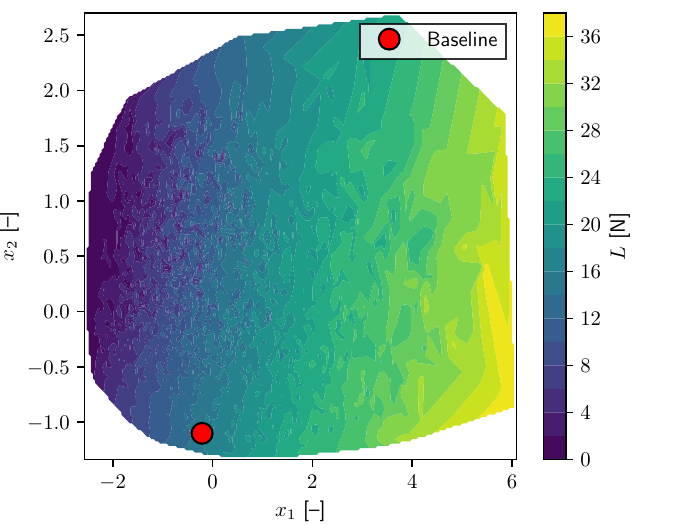}
    }
    \subfigure[TC\#3: Propeller]{
    \includegraphics[width=0.31\linewidth]{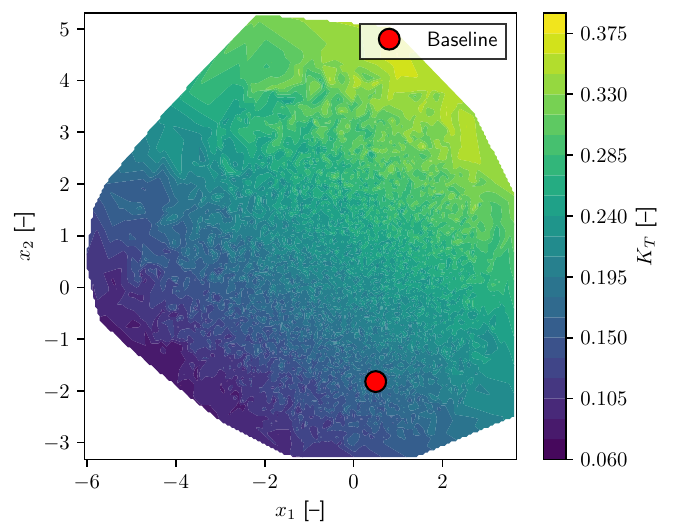}
    \includegraphics[width=0.31\linewidth]{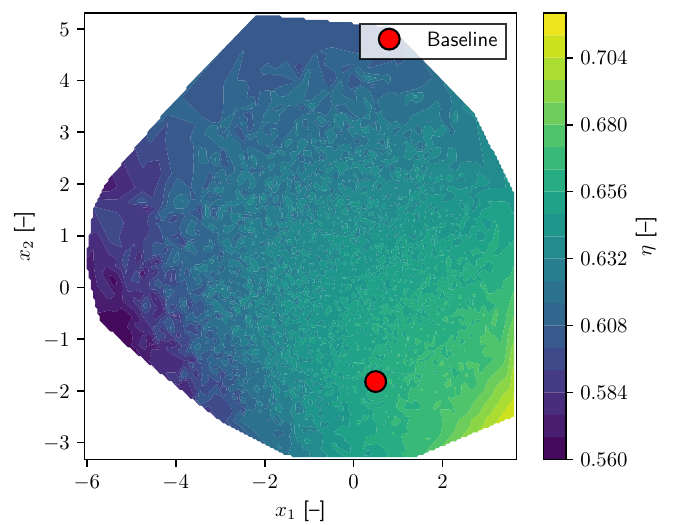}
    \includegraphics[width=0.31\linewidth]{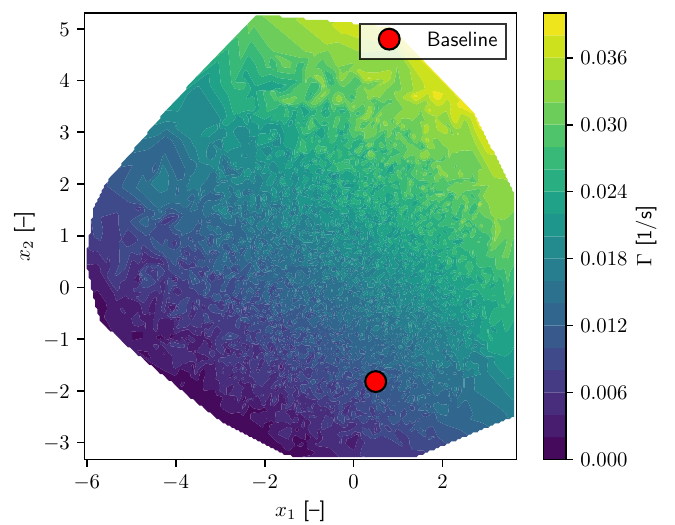}}
    \subfigure[TC\#4: 5415]{
    \includegraphics[width=0.31\linewidth]{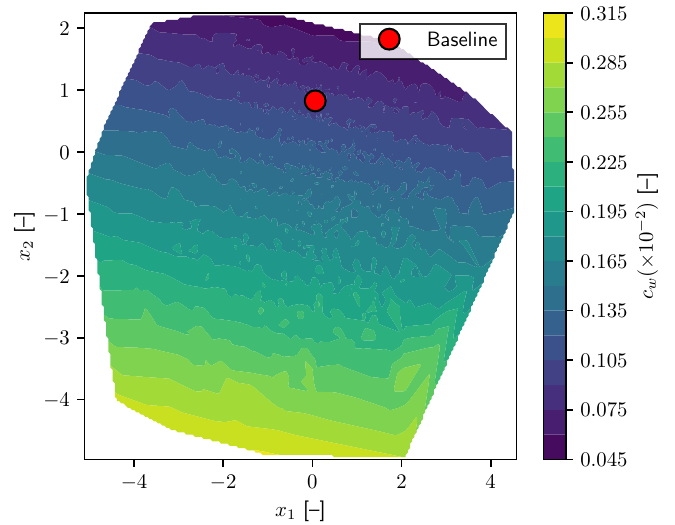}
    \includegraphics[width=0.31\linewidth]{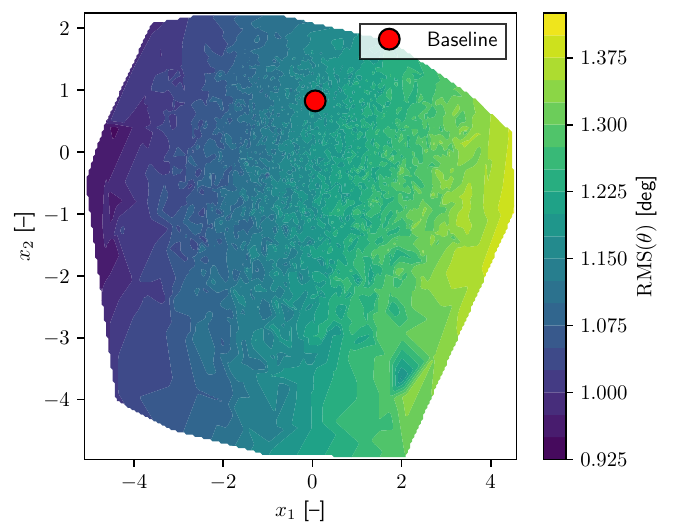}
    \includegraphics[width=0.31\linewidth]{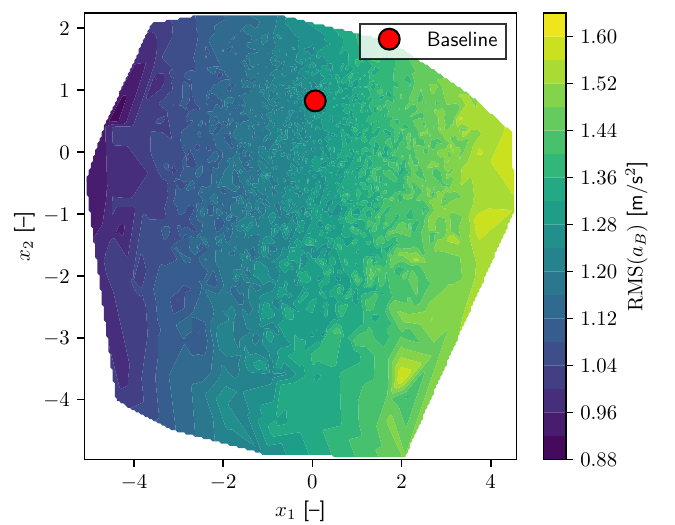}}
    \caption{Contour maps of lumped physical information as a function of the first two reduced design variables $\bfx$}
    \label{fig:contours}
\end{figure}
\begin{figure}[!b]
    \centering
    \subfigure[RAE-2822]{\includegraphics[width=0.44\linewidth]{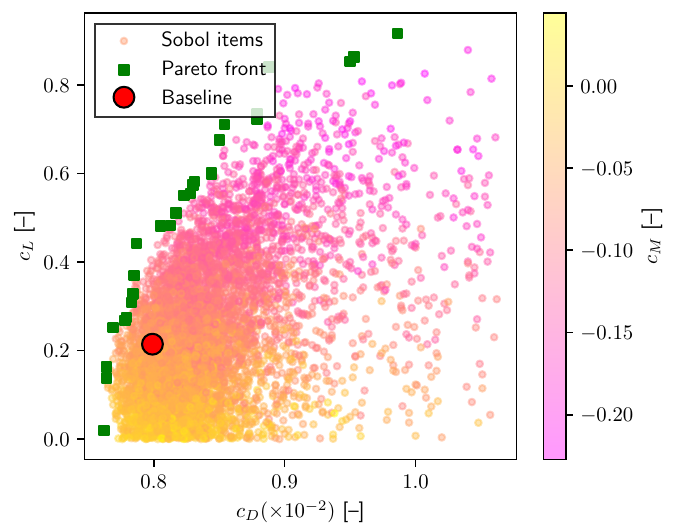}}
    \subfigure[AUG]{\includegraphics[width=0.44\linewidth]{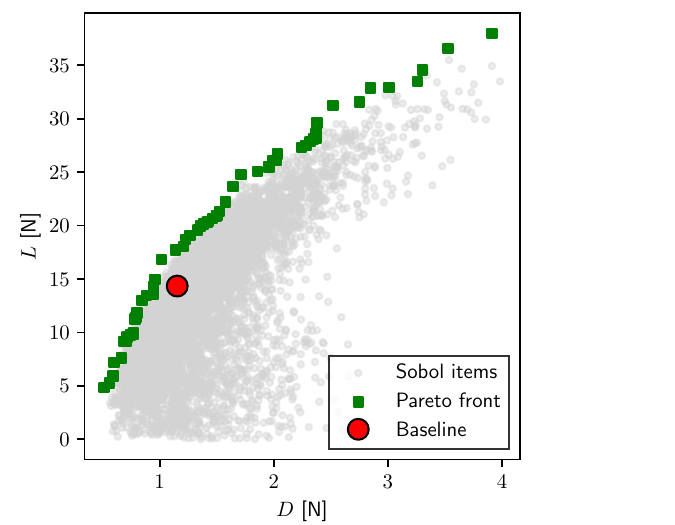}}\\
    \subfigure[Propeller]{\includegraphics[width=0.44\linewidth]{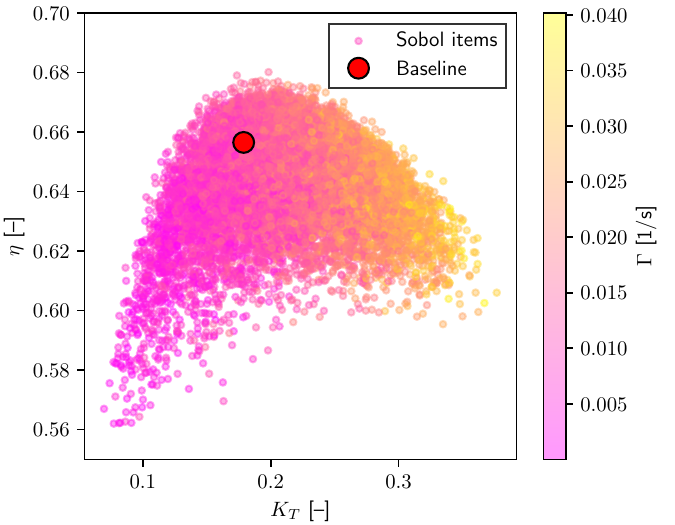}}
    \subfigure[DTMB-5415]{\includegraphics[width=0.44\linewidth]{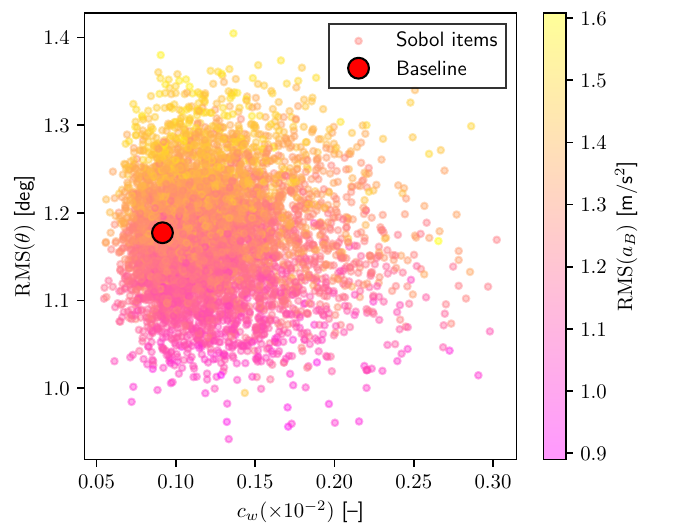}}
    \caption{Scatter plot of lumped physical parameters usable as objective functions in a multi-objective problem}
    \label{fig:paretos}
\end{figure}
\section{Discussion}
\change{
One of the primary motivations behind the proposed PI-/PD-PME extensions is to incorporate physical insight into the dimensionality reduction process without imposing prohibitive computational costs. While enriching the embedding with physical quantities (e.g., lift, drag, or pressure fields) might suggest an increased overhead, these quantities are obtained from low-fidelity solvers (e.g., \texttt{XFOIL} and \texttt{PUFFIn}), which are computationally inexpensive and widely used in early-stage design. Importantly, in many practical workflows, these evaluations are already available from previous parametric studies, and thus no new high-fidelity simulations are required for training. This enables a cost-effective embedding that reflects relevant physical trends.
}

\change{
This rationale is further supported by an empirical analysis of variance convergence. For all test cases considered, the evolution of both geometric and physical variance was monitored as a function of the number of Sobol samples used in the training dataset. Despite generating up to 16,384 samples, it was observed that variance saturation typically occurs around 1,000 samples. This indicates that the main informative directions in the design space can be captured with a limited number of low-fidelity simulations. Such sampling budgets are compatible with standard preliminary design campaigns, even for problems involving more than 20 design variables. The convergence behavior is illustrated in Fig.~\ref{fig:var_conv}, where both geometric and physical variance show rapid stabilization across all configurations. Once this saturation is reached, the structure of the reduced manifold becomes largely insensitive to additional samples, supporting the robustness of the method.}

\change{Moreover, investing in a physics-informed latent space leads to significant downstream savings. By embedding performance-related directions early, the method avoids wasting computational effort on uninformative or non-promising shape perturbations.
This is confirmed, for instance, in the propeller case study, where the dominant reduced basis direction consistently aligns with pitch at tip, a shape feature known to have strong influence on hydrodynamic performance. Similar trends are observed across other configurations, indicating that the reduced space prioritizes design variations with high physical relevance.
This is particularly valuable when transitioning to high-fidelity or multi-objective optimization, where the cost of full evaluations is orders of magnitude higher. Unlike purely geometric PCA or POD, which may generate reduced spaces uncorrelated with performance, or black-box regressors that do not provide a true reduction of the design space, PI-/PD-PME yield interpretable, physics-aware coordinates that guide the design exploration more effectively.}

\change{It is also important to note that the proposed dimensionality reduction does not require highly accurate physical outputs: it relies on capturing global variance trends rather than absolute values. As long as the low-fidelity solver provides coherent gradients and relative sensitivities across the design space, the method effectively identifies the most informative directions. 
This assumption is supported by common practice in design optimization, where low-fidelity models are often benchmarked against high-fidelity responses along selected design directions to validate their qualitative behavior. Moreover, the dimensionality of the reduced space is deliberately kept higher than strictly necessary to mitigate potential inaccuracies in the low-fidelity-based embedding and to retain sufficient flexibility for high-fidelity refinement. 
This allows the embedding to serve not only as a compact representation of the space, but also as a candidate low-fidelity model in a multi-fidelity optimization framework.}
Furthermore, early sensitivity analyses are made easier by looking at, e.g., the 2D contour plots (see Fig. \ref{fig:contours}), which reveal how varying just one or two reduced coordinates shifts integral performance (e.g., drag, thrust, seakeeping motion) away from the baseline configuration.
This ``lumped-level'' lens is invaluable for quick screening, where high-fidelity CFD simulations are too expensive to run exhaustively. By building on these low-fidelity models, designers can prune unpromising shape variations early on. Furthermore, this low-fidelity model can be already used to illustrate the trade-offs between multiple objectives (e.g., lowering drag vs. improving lift, increasing thrust vs. reducing vortex intensity, or minimizing calm-water resistance while controlling seakeeping motions), as well as providing a first Pareto front of optimal solutions (see Fig. \ref{fig:paretos}).

This approach allows for quickly visualizing constraints, identifying high-potential regions, and serving as a stepping stone, guiding more elaborate, high-dimensional or high-fidelity optimization, focusing computational resources on the most promising design clusters.

Alongside lumped metrics, the embedding often includes distributed fields—pressure over the surface for the airfoil and glider, or near-blade flows for the propeller. These data prove critical if the final goal is to train machine-learning reduced-order models (ROMs) capable of predicting integral and/or local quantities of interest. For instance,  the pressure distribution, for the RAE-2822, can be aggregated into a small set of modes that reconstruct aerodynamic loads, while, for the propeller, it provides richer mode shapes and ensures variations that affect vortex shedding, cavitation risk, or tip flow details. %are captured.
Since the embeddings reduce the dimension of these distributed fields, surrogate ROMs can be trained more efficiently, predicting crucial integral outputs without re-running the solver for every shape perturbation.

\change{
Finally, a distinctive feature of the PI-/PD-PME framework is its ability to accommodate a large number of physical quantities without requiring prior manual selection. When the most informative observables are not known in advance, the method identifies dominant contributors based on their empirical variance and statistical correlations. It should be noted, however, that physical metrics that are redundant or irrelevant for performance are not automatically downweighted: if repeated or exhibiting high variance, they may still dominate the embedding. As such, proper normalization and careful feature selection remain important to avoid unintentional overemphasis. Nonetheless, this unsupervised, data-driven formulation enables the identification of performance-relevant directions and supports the exploration of complex design spaces where physical intuition alone may be insufficient. The relative impact of each physical feature on the reduced modes is visualized in Figs.~\ref{fig:rae_eigs}–\ref{fig:5415_eigs}.
}

\section{Conclusions}
This work has demonstrated how the parametric model embedding (PME) framework can be extended to include physical information—either partially (PI-PME) or exclusively (PD-PME)—in a manner that seamlessly combines geometry and performance data within a reduced design space for shape optimization. Across a range of test cases, including the RAE-2822 airfoil, a bio-inspired underwater glider, a ship propeller, and the DTMB-5415 destroyer-type hull, the results showed that adding physics to PME does not always raise the number of modes needed for a target variance level, especially when physical behavior strongly correlates with geometric modifications. In cases where geometry and physics are poorly aligned, however, PI-PME incorporates both purely geometric and purely physical directions, inevitably increasing the dimensionality. This difference reveals how important it is to identify whether the chosen geometric parameterization meaningfully covers the performance-critical aspects of the design.

PI-PME preserves the comprehensive interplay between shape and physics and is particularly useful when geometric feasibility, manufacturing requirements, and interlinked performance metrics must remain in the same design space. PD-PME, by removing geometry, directly focuses on physics-driven changes and can yield a reduced representation when many shape variations do not meaningfully affect integral objectives such as resistance, lift, or thrust. In this role, PD-PME also serves as a diagnostic to highlight situations where geometric degrees of freedom provide little to no performance advantage. Lumped physical quantities—such as lift, drag, thrust, or calm-water resistance—can be incorporated at lower fidelity for early screening or coarse optimization, while distributed fields like surface pressure or vortex data enrich the understanding of local flow phenomena and facilitate training machine-learning-based reduced-order models for predicting global performance measures.

A clear benefit of these methods lies in their interpretability. By revealing which original design variables contribute most strongly to each embedding mode, the extended PME approaches guide engineers to focus on shape changes that matter for performance, rather than for mere geometric variety. Early-stage screening becomes more efficient, as shown by low-dimensional contour plots and Pareto fronts generated from PI-PME or PD-PME analyses, allowing for fast identification of promising design regions without incurring the costs of exhaustive high-fidelity simulation. Ultimately, whether one relies on PME alone (when geometry largely dictates performance), integrates physical data using PI-PME (when shape–physics coupling is important), or centers on PD-PME (when physics alone can drive meaningful decisions), these variations collectively form a flexible toolkit for modeling, analyzing, and optimizing complex design spaces.

Finally, while the proposed methodology effectively incorporates physical information into the reduced space, a limitation lies in its inherently linear structure. In practical applications, both geometric deformations and physical responses often exhibit nonlinear behavior that may not be fully captured by linear embeddings.
PME does not impose explicit assumptions regarding the smoothness or continuity of the geometric or physical fields. Its basis construction is entirely data-driven and reflects the empirical variance across samples. As a linear technique, PME can still capture sharp gradients or localized features—provided they manifest as dominant directions of variability. However, in design spaces characterized by strong discontinuities or highly nonlinear relationships, its effectiveness may be limited, and nonlinear extensions could provide more expressive representations.
Therefore, future research should explore nonlinear extensions of PME by leveraging advanced dimensionality reduction techniques, such as local PCA \cite{d2018nonlinear}, kernel PCA \cite{gaudrie2020modeling,zhao2024supervised}, autoencoders \cite{seo2024study,abbas2023deep,yamazaki2023comparative,karafi2024simultaneous,boncoraglio2021active,kou2023aeroacoustic}, or principal geodesic analysis \cite{bo2024data,doronina2025aerodynamic}, which allow for the discovery of curved low-dimensional manifolds. These techniques may further enhance the ability of the framework to capture complex dependencies between design variables and physical observables, improving both the fidelity and generalizability of the reduced representations.
Nonetheless, to the best of the authors' knowledge, none of these nonlinear approaches currently support the core capabilities offered by PME—namely, the ability to assign custom weights across heterogeneous features and to perform analytical backmapping to the original parametric space. Developing nonlinear extensions of PI-/PD-PME that retain such properties represents a compelling and largely unexplored direction for future research.

%%%%%%%%%%%%%%%%%%%%%%%%%%
\section*{Acknowledgements}
This work has been conducted within the NATO-AVT-404 Research Task Group on ``Enhanced Design Processes of Military Vehicles through Machine Learning Methods''. CNR-INM and the University of Genoa acknowledge the support of the Italian Ministry of University and Research (MUR) through the National Recovery and Resilience Plan (PNRR), Sustainable Mobility Center (CNMS), Spoke 3 Waterways, CN00000023 - CUP B43C22000440001. CNR-INM is also grateful to MUR through PRIN 2022 program, project BIODRONES, 20227JNM52 - CUP B53D23005560006. CNR-INM finally acknowledges the CINECA award under the ISCRA initiative, for the availability of high-performance computing resources and support. CIRA work was supported by the CIRA internal project OPTIWING (OPTImization for WINg Generation).

%Bibliography
\bibliographystyle{unsrt}  
\bibliography{biblio}

\end{document}